\documentclass[fleqn,10pt]{wlscirep}

\usepackage{color,array}

\usepackage{graphicx}
\usepackage{algorithm,algorithmic}
\usepackage{tabularx}
\newcolumntype{Y}{>{\centering\arraybackslash}X}
\usepackage{multirow}
\usepackage{amsmath}
\usepackage{amssymb}
\usepackage{xcolor}
\usepackage[normalem]{ulem}
\usepackage{subcaption}

\usepackage{diagbox}
\usepackage{makecell}
\usepackage{cite}

\newcommand{\mmh}[0]{mm$^2\,$h$^{-1}$}

\newcommand{\commentout}[1]{}

\newcommand{\supplref}[2]{#1 \ref{#2} in the appendix}

\usepackage[nameinlink,noabbrev,capitalize]{cleveref}
\usepackage{subcaption}

\usepackage{xr-hyper}

\makeatletter
\newcommand*{\addFileDependency}[1]{
  \typeout{(#1)}
  \@addtofilelist{#1}
  \IfFileExists{#1}{}{\typeout{No file #1.}}
}
\makeatother
\newcommand*{\myexternaldocument}[1]{%
    \externaldocument{#1}%
    \addFileDependency{#1.tex}%
    \addFileDependency{#1.aux}%
}

\myexternaldocument{supplementary}

\DeclareMathAlphabet{\mathcal}{OMS}{cmsy}{m}{n}

\title{Investigating molecular transport in the human brain from MRI with physics-informed neural networks
}
\author[1]{Bastian Zapf}
\author[2]{Johannes Haubner}
\author[2]{Miroslav Kuchta}
\author[3,4]{Geir Ringstad}
\author[5,6]{Per Kristian Eide}
\author[1,2, *]{Kent-Andre Mardal}
\affil[1]{University of Oslo, Faculty of Mathematics and Natural Sciences, Oslo, 0851, Norway}
\affil[2]{Simula Research Laboratory, Department of Numerical Analysis and Scientific Computing, Oslo, 0164, Norway}

\affil[3]{Oslo University Hospital, Department of Radiology, Oslo, 0372, Norway}
\affil[4]{Sorlandet Hospital, Department of Geriatrics and Internal medicine, Arendal, 4838, Norway}

\affil[5]{Oslo University Hospital, Department of Neurosurgery, Oslo, 0372, Norway}
\affil[6]{University of Oslo, Institute of Clinical Medicine, Oslo, 0372, Norway}

\affil[*]{kent-and@simula.no}

\begin{abstract}
In recent years, a plethora of methods combining neural networks and partial differential equations have been developed. 
A widely known example are physics-informed neural networks, which solve problems involving partial differential equations by training a neural network.
We apply physics-informed neural networks and the finite element method to estimate the diffusion coefficient governing the long term spread of molecules in the human brain from magnetic resonance images.
Synthetic testcases are created to demonstrate that the standard formulation of the physics-informed neural network faces challenges with noisy measurements in our application.
Our numerical results demonstrate that the residual of the partial differential equation after training needs to be small for accurate parameter recovery.
To achieve this, we tune the weights and the norms used in the loss function and use residual based adaptive refinement of training points.  
We find that the diffusion coefficient estimated from magnetic resonance images with PINNs becomes consistent with results from a finite element based approach when the residuum after training becomes small.
The observations presented here are an important first step towards solving inverse problems on cohorts of patients in a semi-automated fashion with physics-informed neural networks.
\end{abstract}

\begin{document}

\flushbottom
\maketitle

\section{Introduction}
In the recent years there has been tremendous 
activity and developments in combining machine learning with physics-based models in the form of partial differential equations (PDE).
This activity has lead to the emergence of the discipline "physics-informed machine learning"  \cite{karniadakis_physics-informed_2021}. 
Therein, nowadays, arguably one of the most popular approaches are physics-informed neural networks (PINNs) \cite{raissi_physics-informed_2019, cuomo2022scientific2}.
They combine PDE and boundary/initial condition into a non-convex optimization problem which can be implemented and solved using mature machine learning frameworks while easily leveraging modern hardware (e.g. GPU-accelerators). 
One of the benefits of the PINN compared to traditional numerical methods for PDE is that no mesh is required. Further, inverse PDE problems are solved in the same fashion as forward problems in PINNs. The only modifications to the code are to add the unknown PDE parameters one seeks to recover to the set of optimization parameters and an additional data-discrepancy term to the objective function.
Among other approaches\cite{rudy_data-driven_2017, peng_multiscale_2021}, PINNs can be used to discover unknown physics from data.
In the context of computational fluid dynamics, PINNs have been successfully applied in inverse problems using simulated data, see, e.g., \cite{cai_heat_2020, jin_nsfnets_2021, reyes_learning_2020, arzani_uncovering_2021} and real data \cite{cai_flow_2021-1, kissas_machine_2020}. A comprehensive review on PINNs for fluid dynamics can be found in  \cite{cai_physics-informed_2021}.

In this work, we solve an inverse biomedical flow problem in 4D with unprocessed, noisy and temporally sparse MRI data on a complex domain. 
Classical approaches require careful meshing of the brain geometry and making assumptions on the boundary conditions \cite{valnes_apparent_2020}.
In patient-specific brain modeling the meshing is particularly challenging and requires careful evaluation of the generated meshes \cite{mardal2022mathematical}. 
Physics-informed neural networks have been applied for the discovery of unknown physics from data without meshing and without regularization \cite{raissi_physics-informed_2019}. 
This makes the PINN method an appealing and promising approach that avoids major challenges in our application and is therefore well worth investigation.
However, PINNs introduce other challenges such as the choice of the network 
architecture, the optimization algorithm and hyperparameter tuning, e.g., 
weight factors in the loss function.
Nevertheless, it is worth to examine how PINNs perform compared to classical algorithms in our application. 

We aim to perform a computational investigation of the glymphatic theory based on and similar to \cite{ray_quantitative_2021, valnes_apparent_2020} with PINNs. We apply them to model the fluid mechanics involved in brain clearance. Various kinds of dementia have recently been linked to a malfunctioning waste-clearance system - the so-called glymphatic system~\cite{iliff_paravascular_2012}. 
In this system, peri-vascular flow of cerebrospinal fluid (CSF) plays a crucial role either through bulk flow, dispersion or even as a mediator of pressure gradients through the interstitium~\cite{nedergaard_glymphatic_2020}. 
While imaging of molecular transport in either rodents~\cite{mestre_flow_2018} or humans~\cite{ringstad_brain-wide_2018} points towards accelerated clearance through the glymphatic system, the detailed mechanisms involved in the system are currently debated~\cite{holter_interstitial_2017, hladky_glymphatic_2022, kedarasetti_arterial_2020, ladron2022perivascular, smith_going_2019, ray_analysis_2019}. 

Our approach builds on previous work where the estimated apparent diffusion coefficient (ADC) for the distribution of gadobutrol tracer molecules over 2 days, as seen in T1-weighted magnetic resonance images (MRI) at certain time points, is compared with the ADC estimated from diffusion tensor images (DTI) \cite{valnes_apparent_2020}.
The ADC of gadobutrol was estimated from the T1-weighted images based on simulations using the finite element method (FEM) for optimal control of the diffusion equation. 
The findings were then compared to estimates of the apparent diffusion coefficient based on DTI.
The latter is a magnetic resonance imaging technique that measures the diffusion tensor of water on short time scales, which in turn can then be used to estimate the diffusion tensor for other molecules, such as gadobutrol \cite{valnes_apparent_2020}. The limited amount of available data prevents from quantifying the uncertainty in the recovered parameters, and makes it a challenging testcase for comparing PINNs and finite element based approaches.

Among other works involving physics-informed neural networks and MRI data\cite{fathi2020super, borges2019physics} several works have previously demonstrated the effectiveness of PINNs in inverse problems related to our problem.
PINNs have been applied to estimate physiological parameters from clinical data using ordinary differential equation models\cite{van2022physics}, but we here consider a PDE model.
Parameter identification problems involving MRI data and PDE have been solved using PINNs \cite{kissas_machine_2020, sarabian2022physics}, but the geometries are reduced to 1-D and hence, taking into account the time dependence of the solution, an effectively two-dimensional problem is solved. Both approaches further involve a data smoothening preprocessing step. 

To the best of our knowledge, this work is the first to estimate physiological parameters from temporally sparse, unsmoothened MRI data in a complex domain using a 4-D PDE model with PINNs. 
We start to verify the PINNs approach on carefully manufactured synthetic data, before working on real data. 
The synthetic testcases reveal challenges that occur for the PINNs due to noise in the data and the sensitivity of the neural network training procedure to different choices of hyperparameters. 
For all of the chosen hyperparameter setting, we evaluate the accuracy of the recovered diffusion coefficient based on the value of the PDE and data loss.
For the synthetic test case, as well as for the real test case, it is required to ensure vanishing PDE loss in order to be consistent with the finite element approach. The question on how this is achieved is addressed by heuristics. 
We investigate using the $\ell^1$-norm instead of $\ell^2$-norm for the PDE loss as an alternative to avoid the overfitting. 
We further discuss how to solve additional challenges that arise when applying the PINNs to real MRI data.
Throughout the paper, we solve the problem with both PINNs and FEM. 

\begin{figure*}[t]
    \centering
    \includegraphics[width=0.95\textwidth]{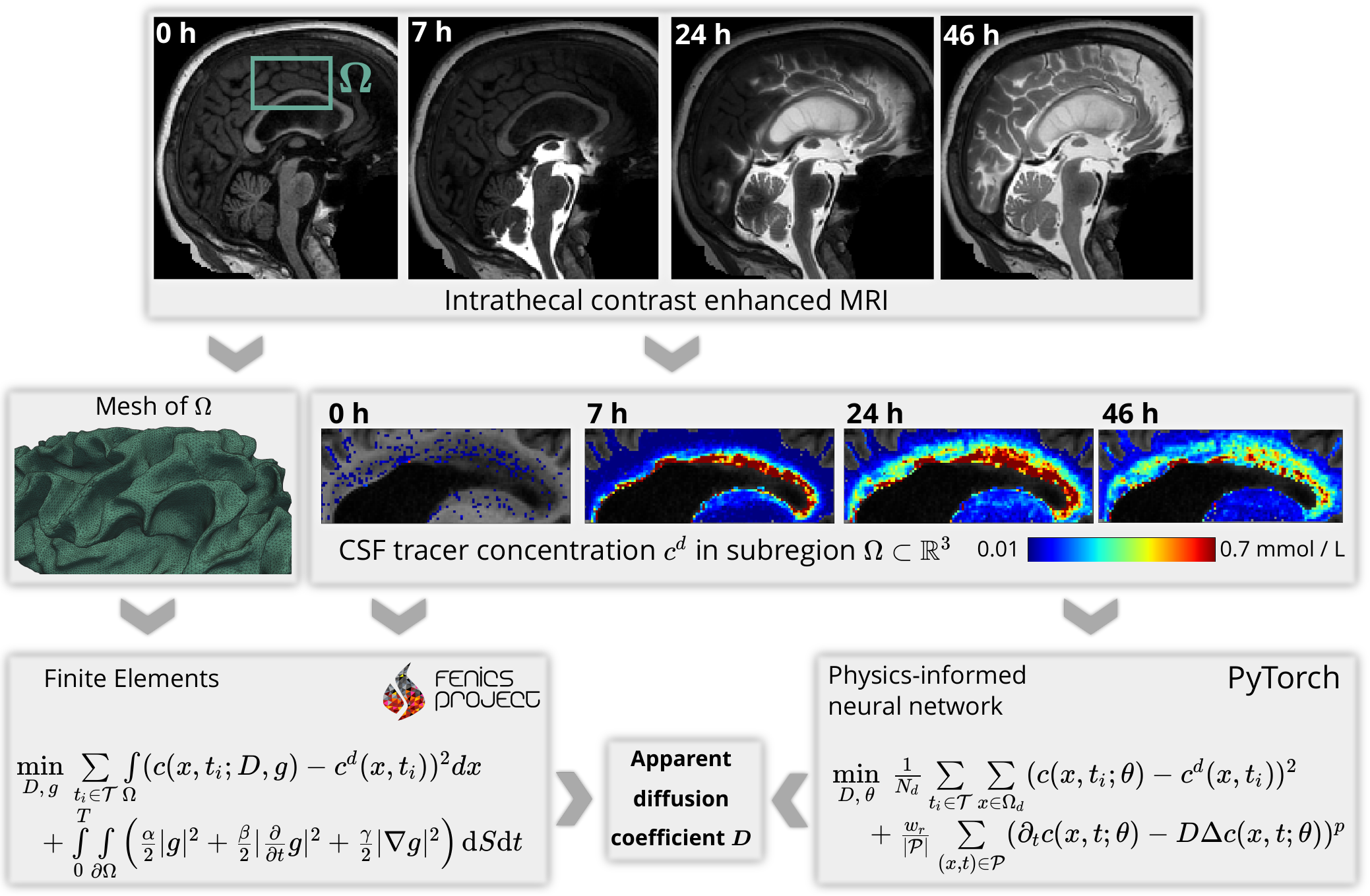}
    \caption{
    Flowchart illustrating our workflow from clinical images to estimated tracer diffusivity in the human brain. From the FreeSurfer \cite{fischl_freesurfer_2012} segmentation of a baseline MRI at $t=0$, we define and mesh a subregion $\Omega$ of the white matter. Intrathecal contrast enhanced MRI at later times $t=7,24, 46$ hours are used to estimate the concentration of the tracer in the subregion. We then use both a finite element based approach and physics-informed neural networks to determine the scalar diffusion coefficient that describes best the concentration dynamics in $\Omega$.
    \label{fig:story}}
\end{figure*}

\section{Problem statement}
Given a set of concentration measurements $c^d(x_j, t_i)$ at four discrete time points 
$t_i \in \{ 0, 7, 24, 46 \}$ hours
and voxel center coordinates $x_j \in \Omega$, where $\Omega \subset \mathbb{R}^3$ represents a subregion of the brain, we seek to find the apparent diffusion coefficient $D>0$ such that a measure
\begin{math}
  J(c, c^d)  
\end{math}
for the discrepancy to the measurement
is minimized under the constraint that $c(x,t)$ fulfills
\begin{align}
    \frac{\partial}{\partial t}  c= D \Delta c  \quad \text{in } \Omega \times (0,T). \label{math:pde} 
\end{align}
The apparent diffusion coefficient takes into account the tortuosity $\lambda$ of the extracellular space of the brain and relates to the free diffusion coefficient $D_f = \lambda^2 D$ \cite{sykova_diffusion_2008}. 
Similar to Valnes \emph{et al.} \cite{valnes_apparent_2020} we here make the simplifying 
assumption of a spatially constant scalar diffusion coefficient.
From the physiological point of view, it is well known that the diffusion tensor in white matter is anisotropic, and hence the modeling assumption of a scalar diffusion coefficient clearly is a simplification. 
The initial and boundary conditions required for the PDE (\ref{math:pde}) to have a unique solution are only partially known, and the differing ways in which we choose to incorporate them into the the PINN and FEM approaches are described in sections \ref{sec:methods-pinn} and \ref{sec:methods-fem}.

\begin{figure}
\centering
     \begin{subfigure}[b]{0.46\textwidth}
         \centering
         \includegraphics[width=\textwidth]{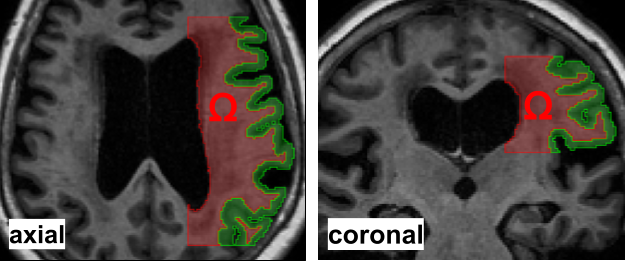}
         \caption{
\label{fig:roi}}
     \end{subfigure}
     \begin{subfigure}[b]{0.46\textwidth}
         \centering
         \includegraphics[width=\textwidth]{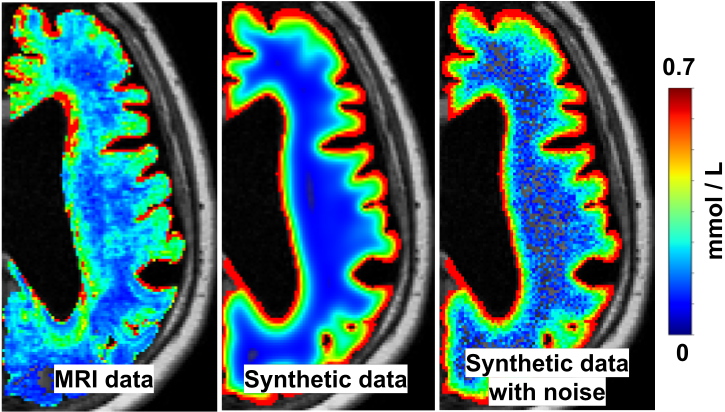}
         \caption{ 
\label{fig:068datasim}}
     \end{subfigure}
    \caption{Geometries and data considered in this work. (a) Axial and coronal slices through the subregion $ \Omega$ of the white matter we consider in this work. The green region depicts the gray matter and is drawn to illustrate the geometrical complexity of the grey matter.
    (b) Axial view of the tracer concentration after 24$\,$h in the right hemisphere for the three data sets considered in this work. Note how the tracer enters the brain from CSF spaces (black). 
    \label{fig:geometry-and-data}
    }
\end{figure}

Our workflow to solve this problem on MRI data is illustrated in Fig. \ref{fig:story}.
Figure \ref{fig:roi} illustrates the white matter subregion $\Omega \subset \mathbb{R}^3$ we consider in this work. 
Figure \ref{fig:068datasim} shows a slice view of the concentration after 24 hours for the three datasets considered in this work, i.e., MRI data, synthetic data with and without noise.
In all cases, we use data at $\mathcal{T} = \{ 0, 7, 24, 46\}$ hours (after tracer injection at $t=0$).

\section{Results}
\subsection{Synthetic data}
We first validate the implementation of both approaches by recovering the known diffusion coefficient $D_0$ from synthetic data without noise. We find that both approaches can be tuned to recover the diffusion coefficient to within a few percent accuracy from three images. Further details can be found in 
\supplref{Section}{suppl-sec:simdata-clean}.

\subsection{Synthetic data with noise \label{results:simdata-w-noise}}

We next discuss how to address challenges that arise for our PINN approach when trained on noisy data as specified by Supplementary Equation (\ref{suppl-math:noise}). 
We find (see \supplref{Table}{suppl-table:minibatch-pinn-clean} for the details) that smaller batch sizes of $\sim\!\!10^4$ points per loss term result in more accurate recovery of the diffusion coefficient (for fixed number of epochs).
We hence divide data and PDE points into 20 batches with $1.5 \times 10^4$ and $5 \times 10^4$ samples per batch, respectively, for the following results.

\begin{figure}
\centering
\captionsetup[subfigure]{labelformat=empty}
     \begin{subfigure}[b]{0.45\textwidth}
         \includegraphics[width=\textwidth]{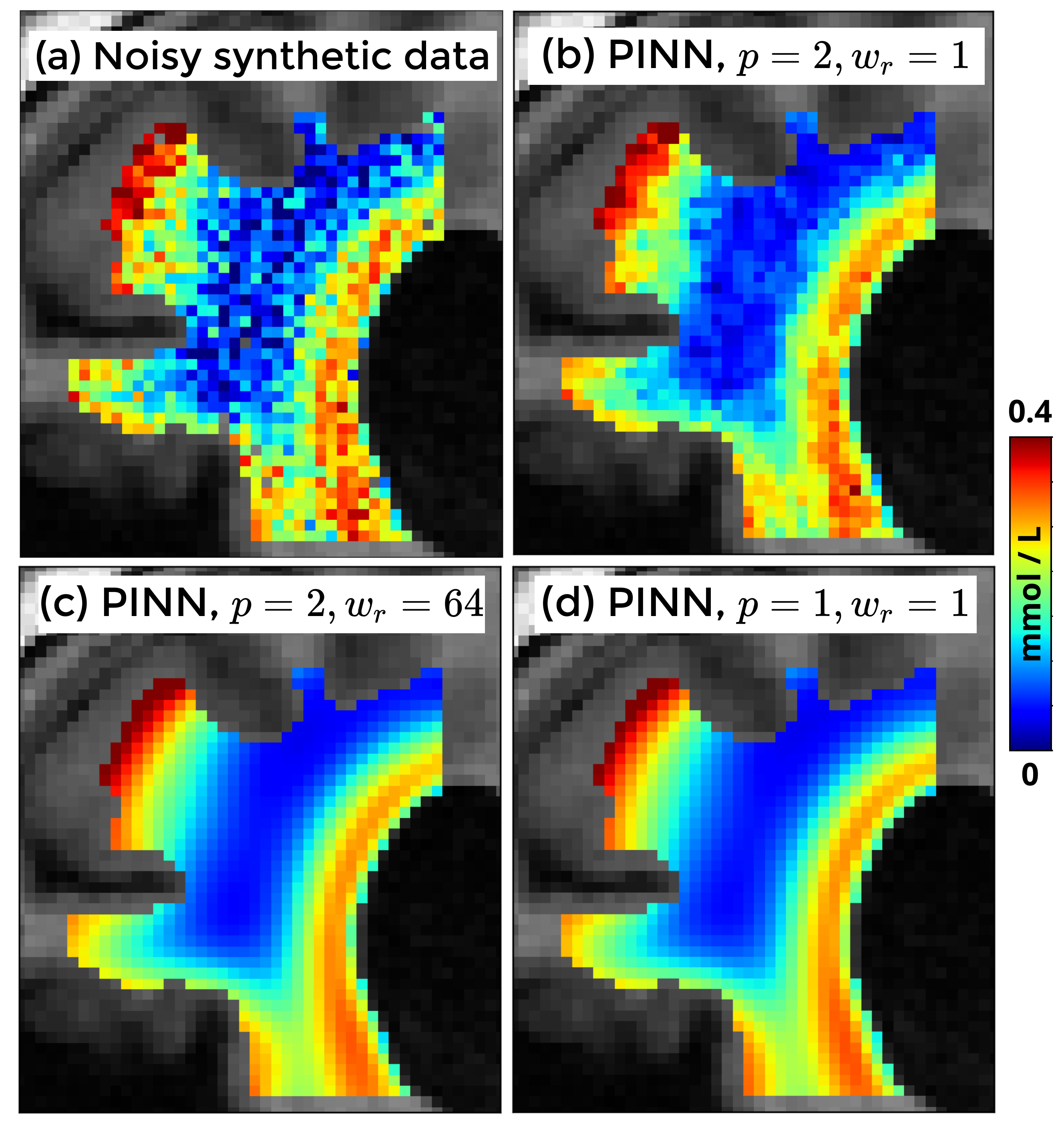}
         \caption{\label{fig:simdata-noise-overfit-a}}
     \end{subfigure}
    \begin{subfigure}[b]{0\textwidth}
         \caption{\label{fig:simdata-noise-overfit-b}}
     \end{subfigure}
         \begin{subfigure}[b]{0\textwidth}
         \caption{\label{fig:simdata-noise-overfit-c}}
     \end{subfigure}
    \begin{subfigure}[b]{0\textwidth}
         \caption{\label{fig:simdata-noise-overfit-d}}
     \end{subfigure}
\captionsetup[subfigure]{labelformat=parens}
\begin{subfigure}[b]{0.46\textwidth}
\vspace{-1cm}
     \begin{subfigure}[b]{\textwidth}
         \centering
         \includegraphics[width=\textwidth]{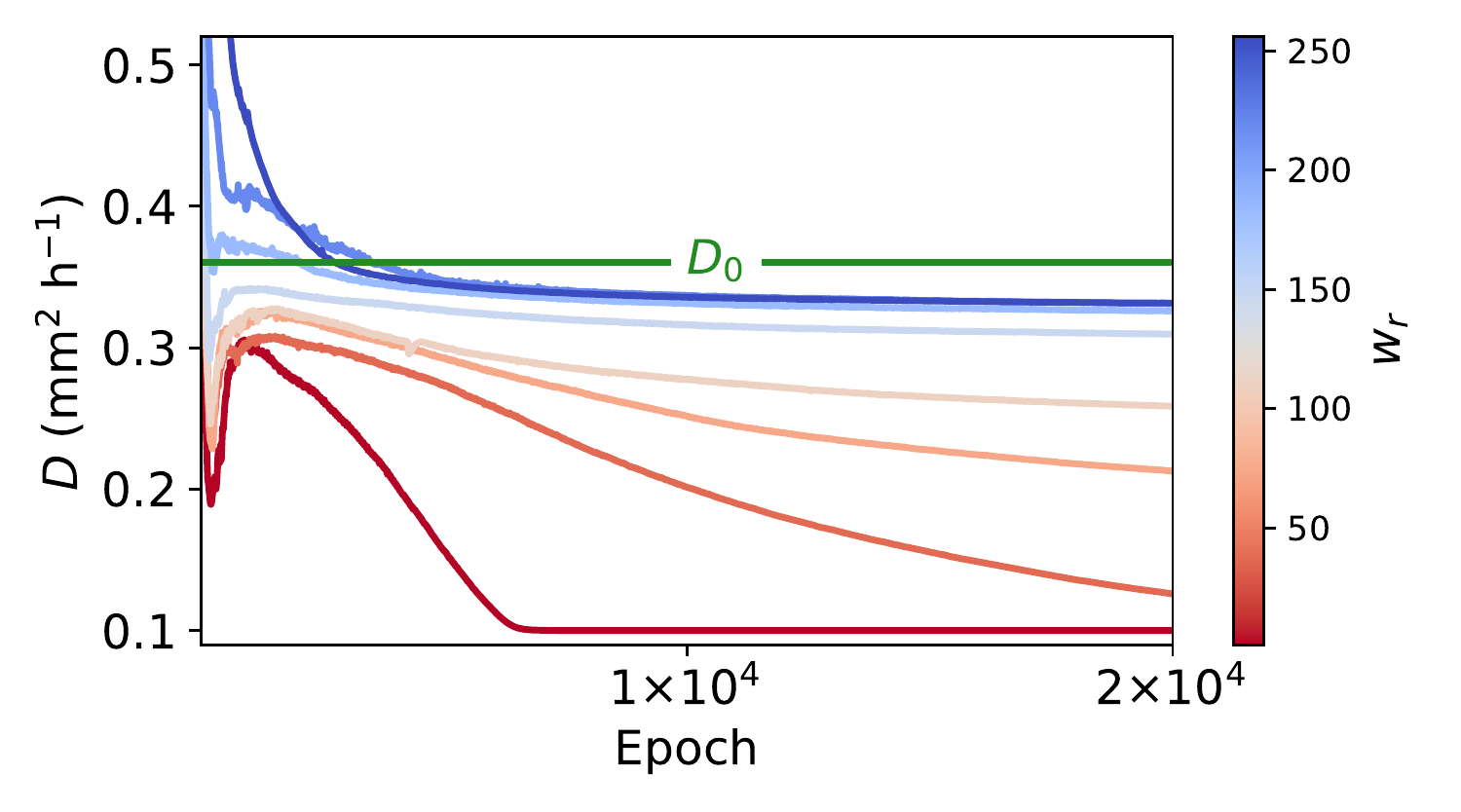}
         \caption{\label{fig:pinn-noise-weights}}
     \end{subfigure}
     \begin{subfigure}[b]{\textwidth}
         \centering
         \includegraphics[width=\textwidth]{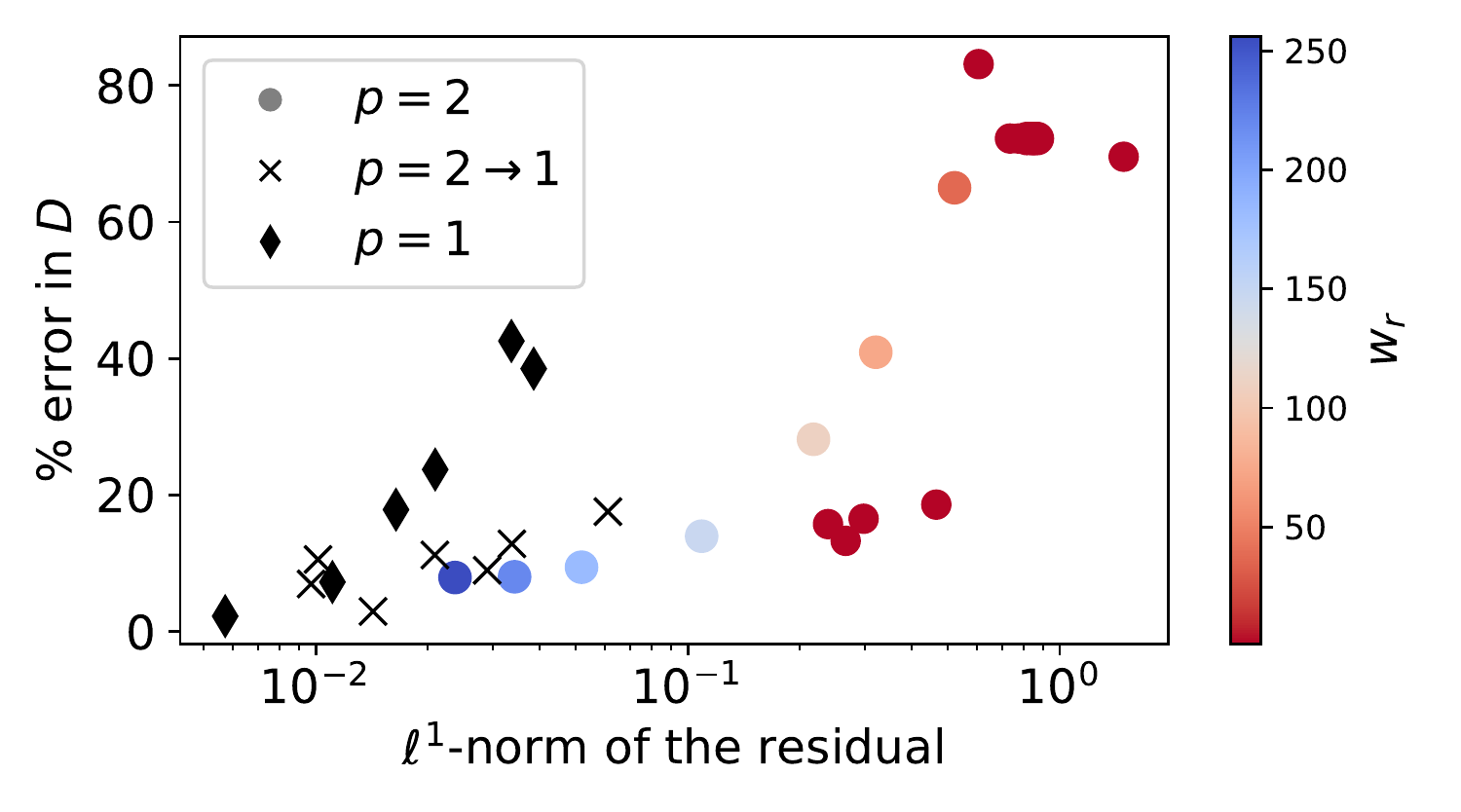}
         \caption{\label{fig:simdata-d-over-res}}
     \end{subfigure}
    \end{subfigure}
    \caption{Influence of PINN hyperparameters on the diffusion coefficient estimated from noisy synthetic data. (a) Coronal slice of synthetic data with noise after $46\,$h, compared to prediction of PINN after training with different hyperparameters in the loss function \eqref{math:pinn}. The overfitting seen in the PINN with $p=2, w_r=1$ (b) can be prevented by using either increased PDE weight $w_r$ (c) or the $\ell^1$-norm for the PDE loss (d).
    (e) 
    The diffusion coefficient recovered by the PINN trained on noisy synthetic data converges to $D_{\mathrm{min}}$ for PDE weight $w_r \leq 2$
    in the loss function \eqref{math:pinn}. 
    (f) Relative error in recovered $D$ from noisy synthetic data as a function of the residual after training for the results presented in (e) and Table \ref{tab:pinn-noise-settings-d}.
    Color encodes the PDE weight $1 \leq w_r \leq 256$ for the results with $p=2$ (dotted). 
    Black markers indicate results with either switching $p=2\rightarrow 1$ during training or $p=1$.
    Different hyperparameter settings in the PINN loss \eqref{math:pinn} yield models which fulfill the PDE to different accuracy, and  low values for the residual coincide with more accurate recovery of the diffusion coefficient.
    \label{fig:noisy-synthetic-data-pinn-results}
}
\end{figure}

In Figs. \ref{fig:noisy-synthetic-data-pinn-results} (a-b) we compare the data to predictions of the PINN after training with the ADAM optimizer\cite{kingma2014adam} and exponential learning rate decay from $10^{-3}$ to $10^{-4}$ for $2\times 10^4$ epochs.
The figures indicate that the network is overfitting the noise that was added to the synthetic data.
This in turn leads to the diffusion coefficient converging to the lower bound $D_{\mathrm{min}}=0.1$ \mmh {} during optimization as shown in Fig. \ref{fig:pinn-noise-weights}.

Here we discuss two remedies:
(i) increasing the regularizing effect of the PDE loss via increasing the PDE weight $w_r$ and
(ii) varying the norm in the PDE loss.
We observe from Fig. \ref{fig:pinn-noise-weights} that for $w_r \gtrsim 64$ the recovered $D$ converges towards the true value to within $\approx 10 \, \%$ error. 
It can also be seen that increasing the weight further does not significantly increase the accuracy.
Fig. \ref{fig:noisy-synthetic-data-pinn-results} (b) and (c) show the predicted solution after 46$\,$h of the trained PINN. 
It can be seen that the overfitting occurring for $w_r=1$ is prevented by choosing a $w_r \geq 64$.
These results are in line with the frequent observation that the weights of the different loss terms in PINNs are critical hyperparameters.
Since we assume that the data is governed by a diffusion equation (with unknown diffusion coefficient), we want the PDE residual to become small.
As demonstrated above, this can be achieved by increasing the PDE weight. 
The correlation between a large weight, a low PDE residual and 
a more accurate recovery of the diffusion coefficient is visualized 
in Fig. \ref{fig:simdata-d-over-res}.

Figure \ref{fig:simdata-d-over-res} also demonstrates the effectiveness the strategy (ii) to successfully lower the PDE residual, which is based on using the $\ell^1$-norm for the PDE loss. 
Using this norm makes the cost function less sensitive to outliers in the data where the observed tracer distribution $c^d$ deviates from the diffusion model \eqref{math:pde}.

Exemplarily, we demonstrate the effectiveness of this approach 
in Fig. \ref{fig:noisy-synthetic-data-pinn-results} (d). There, we plot the PINN prediction after training with $p=1$. It can be seen that the prediction is visually identical to the prediction obtained with $p=2$ and $w_r=64$ (The relative difference between the predictions in Fig. \ref{fig:noisy-synthetic-data-pinn-results} (c) and (d) is about 2$\,$\%). 

The results in Fig. \ref{fig:simdata-d-over-res} are obtained in a systematic study with fixed $w_r=1$. In detail, we test the combinations of the following hyperparameters:
\begin{itemize}
    \item Parameterizations $D(\delta)$ (\ref{math:d-bounded}) vs. $D = \delta$ (\ref{math:d-delta}) of the diffusion coefficient in terms of a trainable parameter $\delta$, c.f. Sec. \ref{sec:dparam}
    \item $p=1$, switching $p=2 \rightarrow 1$ after half the epochs, $p=2$
    \item fixed learning rate $10^{-3}$, exponential learning rate  decay $10^{-3} \rightarrow 10^{-4}$, fixed learning rate $10^{-4}$ and exponential learning rate decay $10^{-4} \rightarrow 10^{-5}$.
\end{itemize}

Table \ref{tab:pinn-noise-settings-d} reports the relative error in the recovered diffusion coefficient after $2 \times 10^4$ epochs of training with ADAM and the minibatch sampling described in \supplref{Algorithm}{suppl-alg:minibatching}.
From the table it can be observed that for $D=\delta$ and $p=1$ instabilities occur with the default learning rate $10^{-3}$ and, due to exploding gradients, the algorithm fails. 
This problem does not occur when using the parameterization $D=D(\delta)$ (\ref{math:d-bounded}). 
It can further be observed that both parameterizations can be fine tuned to achieve errors $\lesssim 10\%$ in the recovered $D$. 
However, the table shows that it is \emph{a priori} not possible to assess which hyperparameter performs best since, for example, settings that fail for the parameterization $D = \delta$ (\ref{math:d-delta}) work well with $D(\delta)$ (\ref{math:d-bounded}).

We hence investigate the effect of the different hyperparameters on the trained PINN and compute
the $\ell^1$-norm of the residual after training defined as
\begin{align}
\frac{1}{{|\mathcal{P}_{\tau}|}} \sum \limits_{(x,t) \in \mathcal{P}_{\tau}} \left| \partial_t c(x,t) -D \Delta c(x,t) \right|.
\label{math:finalresidual-norm}
\end{align}
Here, {$\mathcal{P}_{\tau} = \tau \times \Omega_p$}, where $\tau = \{0, \dots, T\}$ are 200 linearly spaced time points between first and final image at $T=46\,$h and $\Omega_r$ denotes the set of center coordinates of all the voxels inside the PDE domain. 
Note that we evaluate (\ref{math:finalresidual-norm}) with the recovered diffusion coefficient, not with the true $D_0$.
Table \ref{tab:pinn-noise-settings-d} also reports this norm for the different hyperparameter settings.
It can be seen that different hyperparameters lead to different norms of the PDE residual. 
Table \ref{tab:pinn-noise-settings-d} 
reveals that low values of the residual correspond to more accurate recovery of the diffusion coefficient. 
These results are plotted together with the results from Fig. \ref{fig:pinn-noise-weights} in Fig. \ref{fig:simdata-d-over-res} where it can be seen that low PDE residual after training correlates with more accurate recovery of the diffusion coefficient.
This underlines our observation that it is important in our setting to train the PINN such that the norm of the PDE residual is small.

\begingroup
\renewcommand{\arraystretch}{1.5}
\setlength{\tabcolsep}{4pt}
\begin{table}[h]
\centering
\caption{Rel. error $|D-D_0|/D_0$ in \% in the diffusion coefficient
and PDE residual norm after training (in brackets)
for different optimization strategies averaged over 4 trainings on synthetic data with noise. 
It can be seen that the accuracy correlates with the PDE residual after training, i.e. the lower the PDE residual, the more accurate the recovered diffusion coefficient. This relation is further illustrated in Fig. \ref{fig:simdata-d-over-res}. Failure of the algorithm is indicated by the symbol "x".
\label{tab:pinn-noise-settings-d}}
\begin{tabularx}{\textwidth}{c *{6}{Y}}
\hline
 Parameterization   &  \diagbox[]{$p$}{lr} & $10^{-3}$ & \makecell{$10^{-3}$  $\rightarrow 10^{-4}$} & $10^{-4}$ & \makecell{$10^{-4}$  $\rightarrow 10^{-5}$}\\
\hline
   & $1$ & x          & x             & 18    (1.6e-02)       & 43      (3.4e-02)        \\
$D=\delta$   & $2\rightarrow 1$ & x          & 7     (9.7e-03)          & 3    (1.4e-02 )        & 13    (3.4e-02)          \\
   & $2$ & 70  (1.5e+00)         & 83    (6.1e-01)          & 16      (2.4e-01)     & 17      (3.0e-01 )        \\
   \hline
 & 1 & 7       (1.1e-02)     & 2    (5.7e-03)           & 24    (2.1e-02)       & 39      (3.9e-02)        \\
$D=D(\delta)$ & $2\rightarrow 1$ & 11      (2.1e-02)     & 11       (1.0e-02 )       & 9  (2.9e-02)          & 18     (6.1e-02)         \\
 & 2 & 72       (7.3e-01)    & 72      (7.7e-01)        & 13   (2.7e-01 )        & 19  (4.7e-01) \\
\hline
\end{tabularx}
\end{table}
\endgroup

Finally, for the FEM approach, \supplref{Table}{suppl-tab:opt-control-w-noise} tabulates the relative error in the recovered diffusion coefficient for solving 
\eqref{math:fem-functional} with regularization parameters spanning several orders of magnitude. The results are in line with the well-established observation that a sophisticated decrease of the noise level and regularization parameters ensures convergence towards a solution \cite{kaltenbacher_iterative_2008}.
In summary, we find that both methods can be tuned to achieve similar accurate recovery of the diffusion coefficient.

\subsection{MRI Data \label{sec:mriresults}}
We proceed to estimate the apparent diffusion coefficient governing the spread of tracer as seen in MRI images. 
It is worth emphasizing here that our modeling assumption of tracer transport via diffusion with a constant diffusion coefficient $D \in \mathbb{R}$ is a simplification, and that we can not expect perfect agreement between model predictions and the MRI data.
Furthermore, closer inspection of the tracer distribution on the boundary in Fig. \ref{fig:068datasim} reveals that, unlike in the synthetic data, the concentration varies along the boundary in the MRI measurements.
Based on these two considerations it is to be expected that challenges with the PINN approach arise that were not present in the previous, synthetic testcases. 
However, our previous observation that smaller PDE residual correlates with more accurate recovery of the diffusion coefficient serves as a guiding principle on how to formulate and minimize the PINN loss function such that the PDE residual becomes small.

Based on the observation that the parameterization $D=D(\delta)$ avoids instabilities during the optimization, we only use this setting in this subsection. The white matter domain $\Omega$ is the same as in the previous section, and we again divide both data and PDE loss into 20 minibatches. We train for $10^5$ epochs using the ADAM optimizer with exponentially decaying learning rate $10^{-4}$ to $10^{-5}$.

\begin{figure}[h]
\centering
     \begin{subfigure}[b]{0.46\textwidth}
         \centering
         \includegraphics[width=\textwidth]{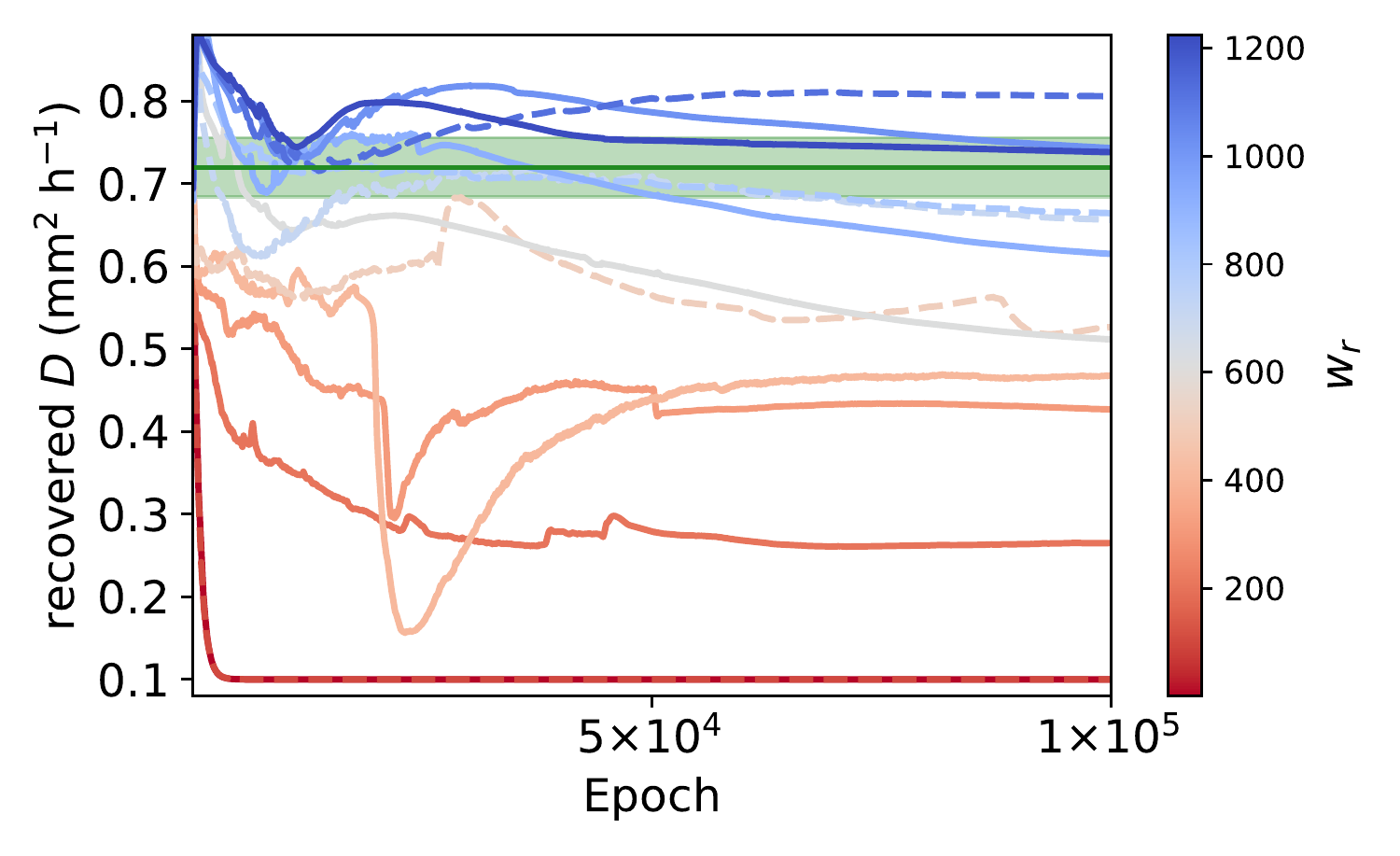}
         \caption{
\label{fig:068_weights}}
     \end{subfigure}
     \begin{subfigure}[b]{0.46\textwidth}
         \centering
         \includegraphics[width=\textwidth]{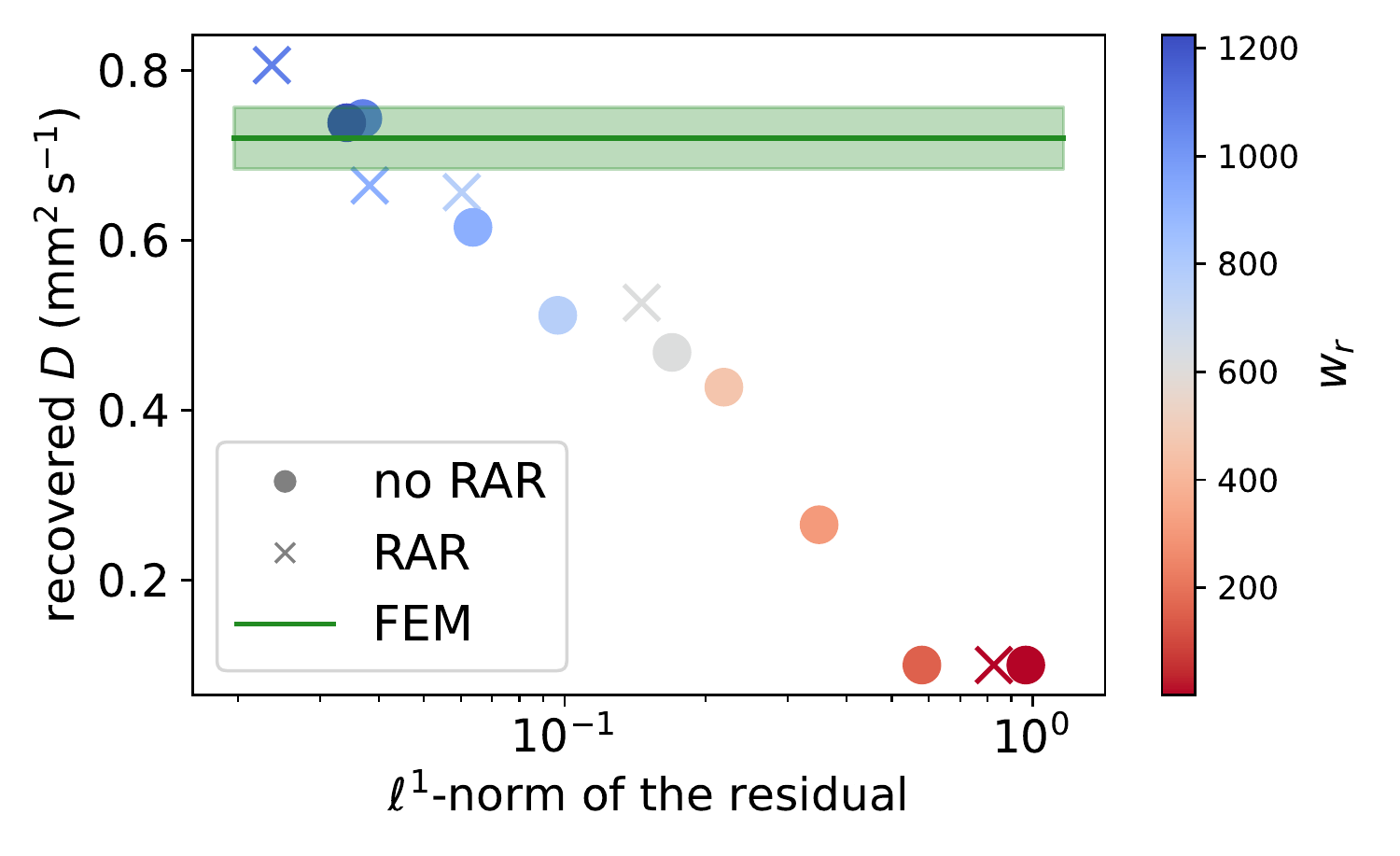}
         \caption{ 
\label{fig:68d-over-residual-weights}}
     \end{subfigure}
    \caption{Influence of PINN hyperparameters on the diffusion coefficient estimated from clinical data.
    (a) 
    Diffusion coefficient during training for different PDE weights $w_r$ and exponentially decaying learning rate from $10^{-4}$ to $10^{-5}$. Dashed lines indicate result with residual based adaptive refinement (RAR).
    (b) Estimated diffusion coefficient with $p=2$ for different PDE weights $w_r$ as a function of the
$\ell^1$-norm of the residual after training. The values for FEM and the green horizontal bars indicating an error estimate are taken from Valnes \emph{et al.} \cite{valnes_apparent_2020}.}
\end{figure}

We first test for $p=2$ with PDE weight 
\begin{math}
    w_r \in \{1, 32, 64, 128, 256, 512, 1024 \}
\end{math}
and display the results in Fig. \ref{fig:068_weights}. It can be seen that, similar to the noisy synthetic data, the diffusion coefficient converges to the lower bound for low PDE weights.
For these settings, we plot the residual norm \eqref{math:finalresidual-norm} of the trained networks in Fig. \ref{fig:68d-over-residual-weights}. It can be seen that increased PDE weight leads to lower residual after training, and in turn to an estimate for $D$ which becomes closer to FEM.

\begin{figure}[h]
\centering
     \begin{subfigure}[b]{0.65\textwidth}
         \centering
         \includegraphics[width=\textwidth]{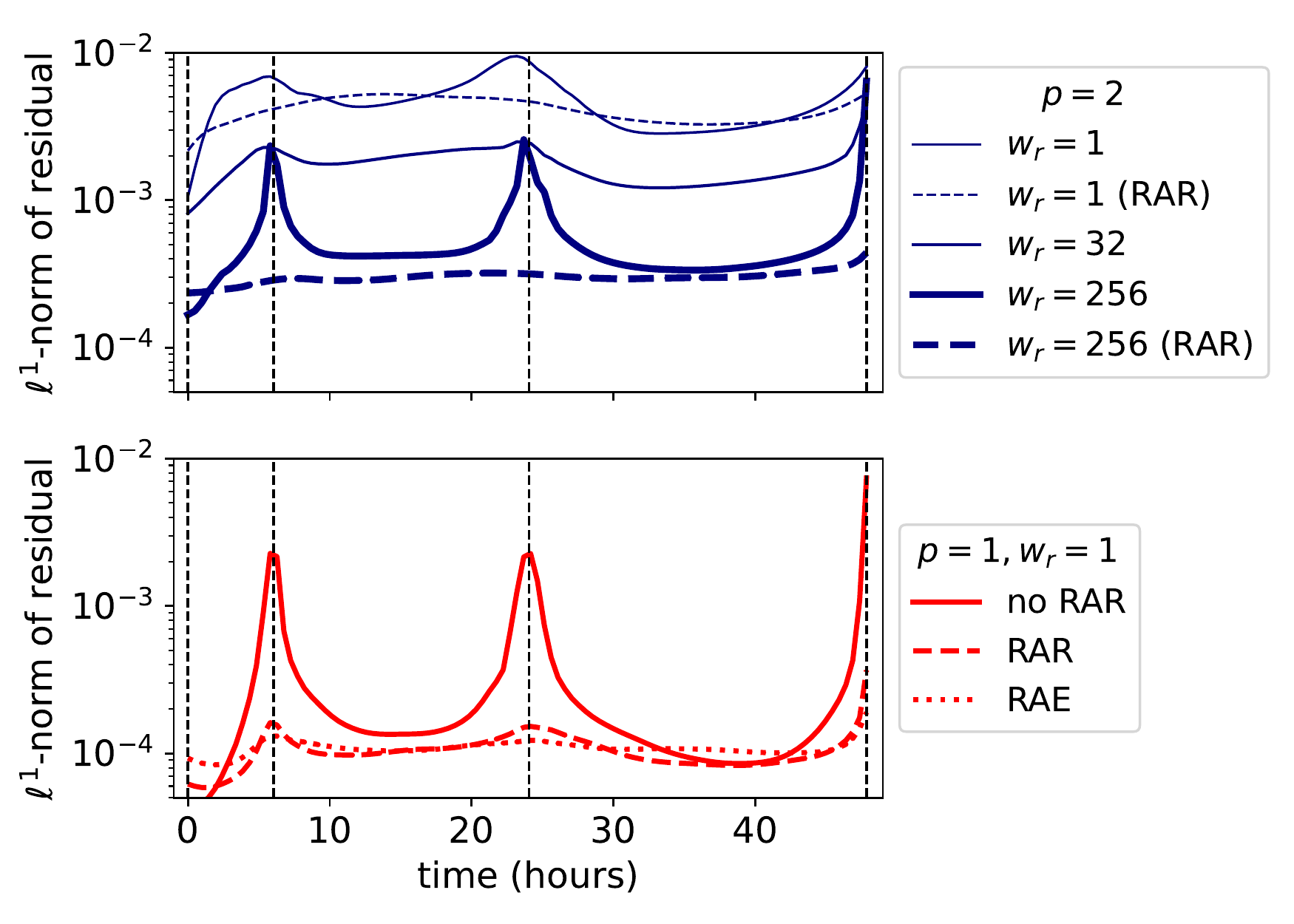}
         \caption{\label{fig:068residual}}
     \end{subfigure}
\begin{subfigure}[b]{0.32\textwidth}
\vspace{-1cm}
\begin{subfigure}[b]{1\textwidth}
         \centering
         \includegraphics[width=0.8\textwidth]{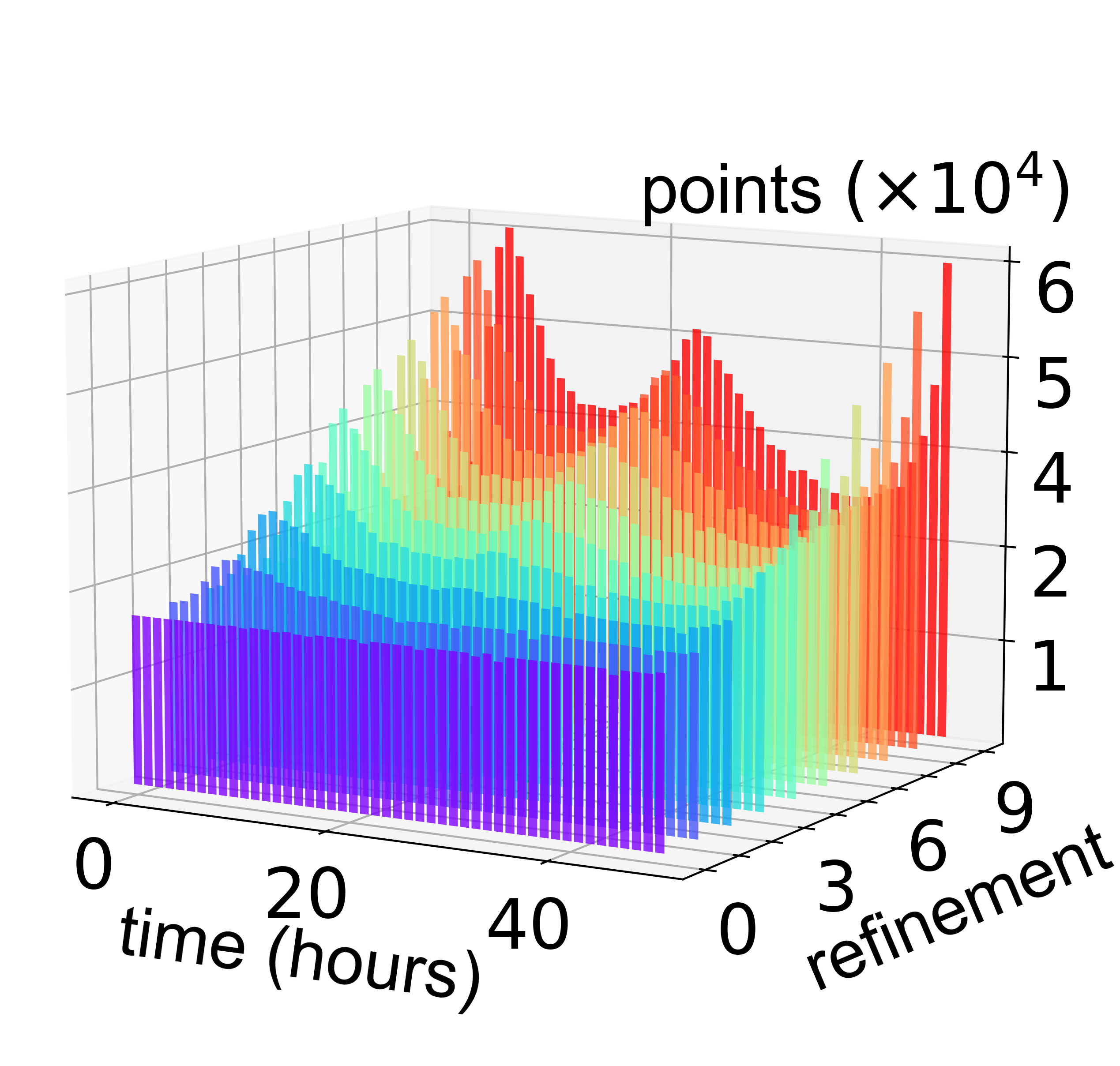}
         \caption{\label{fig:068rar}}
\end{subfigure}
\begin{subfigure}[b]{1\textwidth}
         \centering
         \includegraphics[width=0.8\textwidth]{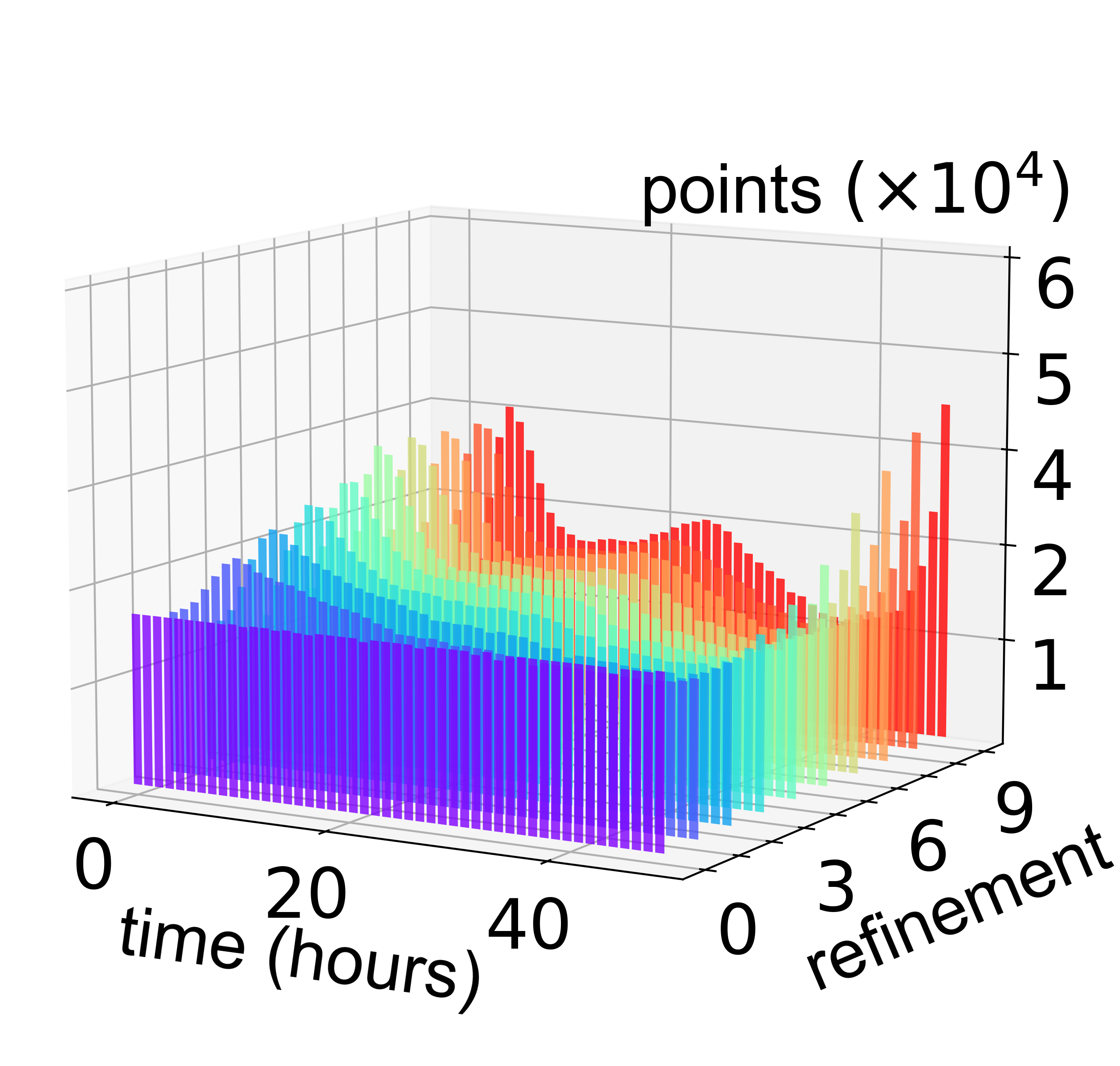}
         \caption{\label{fig:068rae}}
\end{subfigure}
\end{subfigure}
    \caption{
    Adaptive training point refinement is needed to fulfill the PDE in all timepoints.
    (a) 
    Average PDE residual in $\Omega_P$ over time for different optimization schemes. Vertical lines (dashed) indicate the times where data is available. In all cases, the learning rate decays exponentially from $10^{-3}$ to $10^{-4}$.
    (b-c) Distribution of PDE points during training with RAR (b) and RAE (c). Starting from a uniform distribution of points (in time), the procedures add more points at 7, 24 and 46 hours where data is available.}
\end{figure}

Further, in Fig. \ref{fig:068residual} we also plot the $\ell^1$-norm of the residual after training as a function of time $t\in [0, T]$, defined as
\begin{align}
r(t) = \frac{1}{|\Omega_r|} \sum \limits_{x \in \Omega_r} \left| \partial_t c(x,t) -D \Delta c(x,t) \right|.
\label{math:residual-norm}
\end{align}
The continuous blue lines in Fig. \ref{fig:068residual} exemplarily show $r(t)$ for some PDE weights. It can be seen that higher PDE weights lead to lower residuals. 
However, for $w_r = 256$ the PDE residual is significantly higher at the times where data is available than in between. 
We did not observe this behavior in the synthetic testcase. 
Since we want the modeling assumption (\ref{math:pde}) to be fulfilled equally in $\Omega \times [0, T]$, we use residual based adaptive refinement (RAR) \cite{lu_deepxde_2021}.
Using the RAR procedure, we add $10^5$ space-time points to the set $\mathcal{P}$ of PDE points after $1\times 10^4, 2\times 10^4, \dots, 9\times 10^4$ epochs. 
Details on our implementation of RAR and an exemplary loss plot during PINN training are given in \supplref{Section}{suppl-app:RAR-RAE}. 
The effectiveness of RAR to reduce this overfitting is indicated by the dashed blue lines in Fig. \ref{fig:068residual}.

\begin{figure}
\centering
     \begin{subfigure}[b]{0.46\textwidth}
         \includegraphics[width=\textwidth]{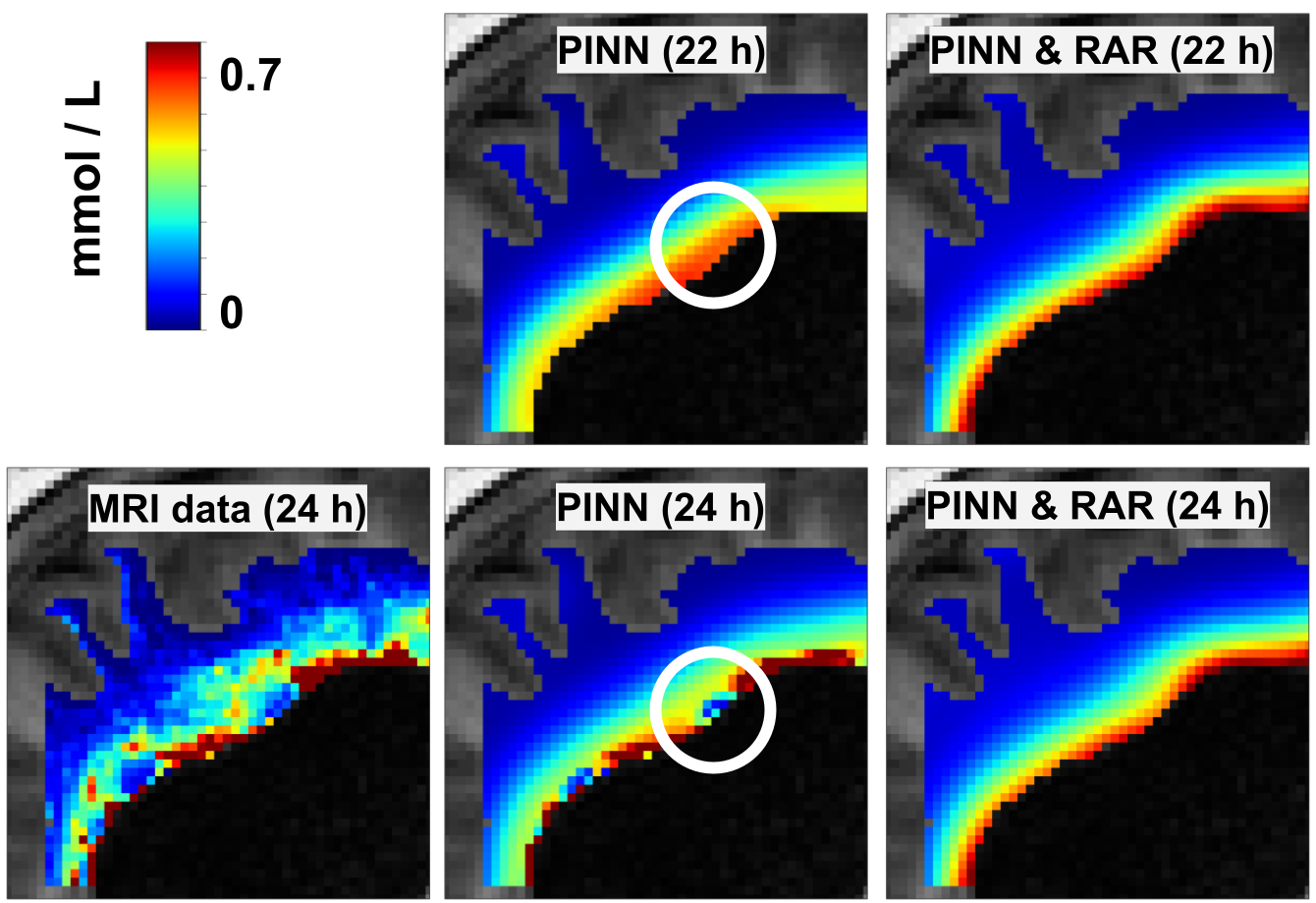}
         \caption{\label{fig:PINN-realdata-overfit}}
     \end{subfigure}
\begin{subfigure}[b]{0.5\textwidth}
         \includegraphics[width=\textwidth]{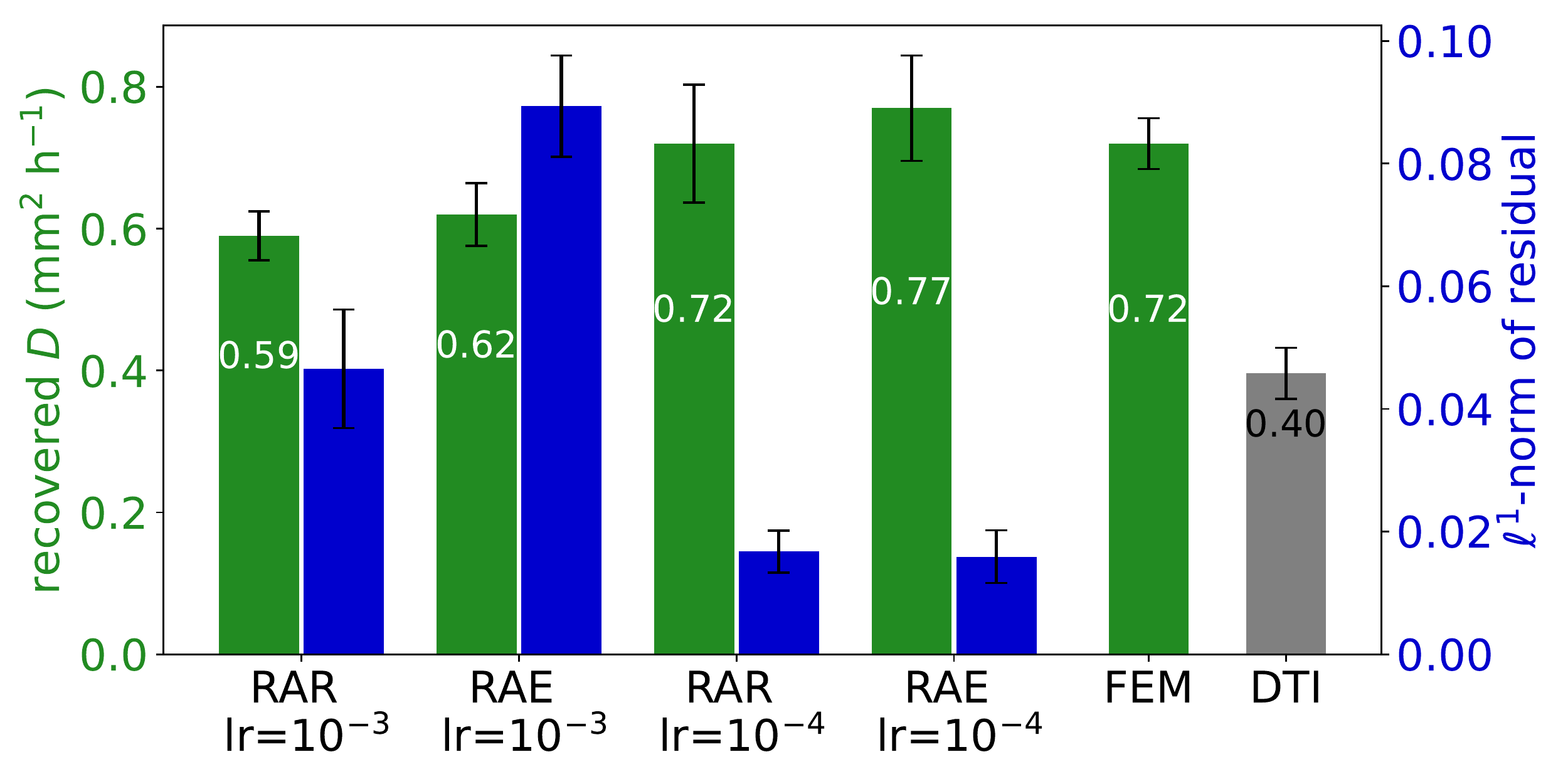}
         \caption{\label{fig:068dhistogram-pinns-vs-valnes}}
\end{subfigure}
    \caption{
    Adaptive refinement yields PINN solutions that are consistent with a diffusion model and FEM results.
    (a)
    Upper row: 
     Prediction of concentration after 22$\,$h from PINNs trained with $p=1$ and $p=1$ \& RAR.
Lower row: Zoom into a sagittal slice of data at $24\,$h compared PINN predictions. 
The PINN prediction after training without RAR overfits the data. Compare also to Fig. \ref{fig:068residual}.
    (b)
    Green: PINN estimates for the diffusion coefficient with RAR or RAE and different initial learning rates ($p=1$ in all cases). 
Blue: $\ell^1$-norm of the residual after training. It can be seen that lower learning rate leads to a lower residual norm and an estimate for the diffusion coefficient closer to the FEM approach. }
\end{figure}

Next, we test for $p=1$ with an exponentially decaying learning rate from $10^{-3}$ to $10^{-4}$ as well as $10^{-4}$ to $10^{-5}$. 
With this setting, the PINNs approach yields an estimate $D= 0.75$ \mmh {} which is close to the FEM solution \cite{valnes_apparent_2020} $D = 0.72$ \mmh {}. 
However, a closer inspecting of the PINN prediction at 22 and 24 (where data is available) shown in Fig. \ref{fig:PINN-realdata-overfit} reveals that the PINN is overfitting the data.
This is further illustrated by the continuous red line in Fig. \ref{fig:068residual} where it can be seen that the PDE residual is one order of magnitude higher at the times where data is available.
The dashed red line in Fig. \ref{fig:068residual} and slices of the predicted $c(x,t)$ shown in Figs. \ref{fig:PINN-realdata-overfit} show that this behavior can be prevented by using RAR.

Since the RAR procedure increases the number of PDE points, the computing time increases (by about 25 \% in our setting). We hence test a modification of the RAR procedure. Instead of only adding points, we also remove the points from $\mathcal{P}$ where the PDE residual is already low. We here call this procedure residual based adaptive exchange (RAE) and give the details in \supplref{Section}{suppl-app:RAR-RAE}.
The dotted red line in Fig. \ref{fig:068residual} demonstrates that in our setting both methods yield similarly low residuals $r(t)$ without overfitting the data. Since in RAE the number of PDE points stays the same during training, the computing time is the same as without RAR. In Figs. \ref{fig:068rar}, \ref{fig:068rae} it can be seen how both RAR and RAE add more PDE points around the timepoints where data is available.

We estimate the apparent diffusion coefficient $D$ by averaging over 5 trainings with either RAR or RAE and learning rate decay from $10^{-3}$ to $10^{-4}$ or $10^{-4}$ to $10^{-5}$.
The results are displayed in Fig. \ref{fig:068dhistogram-pinns-vs-valnes} together with the $\ell^1$-norm (\ref{math:finalresidual-norm}) after training. It can be seen that for the same learning rate, both RAR and RAE yield similar results. 
A lower learning rate, however, leads to lower PDE residual and an estimated diffusion coefficient which is closer to the value 0.72 \mmh {} from Valnes \emph{et al.}  \cite{valnes_apparent_2020}.

\subsection{Testing different patients}
In Valnes \emph{et al.} \cite{valnes_apparent_2020}, the same methodology was applied to two more patients, named 'REF' and 'NPH2'. 
We here test how well the optimal hyperparameter settings found in Section \ref{sec:mriresults} generalize to these patients.
A similar subregion of the white matter is used but the voxels on the boundary of the domain were removed.

A PINN is trained with the following hyperparameters from Section \ref{sec:mriresults} that yielded the lowest PDE residual after training:
The number of minibatches is set to 20, training for $10^5$ epochs with ADAM and exponential learning rate decay from $10^{-4}$ to $10^{-5}$, and $p=1$ with RAR at $1\times 10^4, 2\times 10^4, \dots, 9\times 10^4$ epochs. The network architecture remains the same.
For patient 'NPH2' we find $D=0.48$ \mmh {} while the FEM approach \cite{valnes_apparent_2020} yields $D=0.50$ \mmh {}.
We find $D=0.41 $ \mmh {} for patient 'REF' while the FEM approach \cite{valnes_apparent_2020} yields $D=0.50$ \mmh {}.

\section{Discussion}
We have tested both PINNs and FEM for assessing the apparent diffusion coefficient in a geometrically complex domain, a subregion of the white matter of the human brain, 
based on a few snapshots of T1-weighted contrast enhanced MR images over the course of 2 days. 
Both methodologies yield similar estimates when properly set up, that is; we find that the ADC is in the range (0.6-0.7) \mmh, depending on the method, whereas the DTI estimate is 0.4 \mmh. As such the conclusion is similar to that of Valnes \emph{et al.} \cite{valnes_apparent_2020}. 
With a proper hyperparameter set-up, PINNs are as accurate as FEM and, given our implementation with GPU acceleration, more efficient than our current FEM implementation.

However, choosing such a set-up, i.e., hyperparameter setting, loss function formulation and training procedure, is still a priory not known and challenging. 
An automated way to find a suitable setting is needed. To this end automated approaches such as AutoML \cite{he2021automl} or Meta learning \cite{psaros2022meta}, could be applied in the future. 
Moreover, theoretical guarantees are required, especially in sensitive human-health related applications.

Our results are in line with the frequent observation that the PDE loss weight is an important hyperparameter. 
Several works have put forth methodologies to choose the weights adaptively during training \cite{wang_when_2022, wang_eigenvector_2020, van_der_meer_optimally_2021, maddu_inverse_2022}, but in practice they have also been chosen via trial-and-error \cite{yin_non-invasive_2021, van2022physics, yazdani_systems_2020}. 
However, in settings with noisy data, it can not be expected that both data loss and PDE loss become zero.
The ratio between PDE loss weight and data loss weight reflects to some degree the amount of trust one has in the data and the physical modeling assumptions, i.e., the PDE.
In this work, we have made the modeling assumption that the data is governed by a diffusion equation, and hence require the PDE to be fulfilled.
This provides a criterion for choosing a Pareto-optimal solution if the PINN loss is considered from a multi-objective perspective \cite{rohrhofer2021pareto}.

From the mathematical point of view, we have sought the solution of a challenging nonlinear ill-posed inverse problem with limited and noisy data in both space and time. 
There can thus be more than one local minimum and the estimated solutions depend on the regularization and/or hyperparameters. 
Here, our main observation is that the diffusion coefficient recovered by PINNs approaches the FEM result when the hyperparameters are chosen to ensure that the PDE residual after training is sufficiently small. 

In general, we think that the current problem serves as a challenging test case and is well suited for comparing PINNs and FEM based methods.
Further, since the finite element approach is well-established and theoretically founded it can serve to benchmark PINNs. 
Our numerical results indicate that the norm of the PDE residual of the trained PINN correlates with the quality of the recovered parameter.
This relates back to the finite element approach where the PDE residual is small since the PDE is explicitly solved.  
In our example, we have found that in particular two methodological choices help to significantly lower the PDE-residual in the PINNs approach: $\ell^1$-penalization of the PDE and adaptive refinement of residual points. 

From the physiological point of view, there are several ways to improve upon our modeling assumption of a spatially constant, scalar diffusion coefficient. 
For instance, an estimate of the local CSF velocity can be obtained by the optimal mass transport technique\cite{koundal_optimal_2020}. 
Another approach could be to learn a spatially varying diffusion coefficient or tensor and compare to diffusion tensor imaging data. 
From an implementational point of view, such methods fit well within our current framework since the PINN formulation is comparably easy to implement and the PDE does not have to be solved explicitly.

\section{Methods}

\subsection{Data acquisition}
The approval for MRI observations was retrieved by the Regional
Committee for Medical and Health Research Ethics (REK) of Health Region South-East, Norway (2015/96) and
the Institutional Review Board of Oslo University Hospital (2015/1868) and the National Medicines Agency
(15/04932-7). The study participants were included after written and oral informed consent.

Details on MRI data acquisition and generation of synthetic data can be found in the \supplref{Section}{suppl-sec:data}.

\subsection{The PINN approach\label{sec:methods-pinn}}
In PINNs, our parameter identification problem 
can be formulated as an unconstrained non-convex optimization problem over the network parameters $\theta$ and the diffusion coefficient $D$ as
\begin{align}
    \min_{\theta, D}  \mathcal{J} + w_r \mathcal{L}_r, \label{math:pinn}
\end{align}
where $w_r > 0$ is a weighting factor.
The data loss $\mathcal{J}$ is defined as
\begin{align}
    \mathcal{J}(c, c^d) = \frac{1}{N_d} \sum_{t_i \in \mathcal{T}} \sum_{x \in \Omega_d} (c(x, t_i) - c^d(x, t_i))^2,
    \label{math:pinn-dataloss}
\end{align}
where $\Omega_d$ is a discrete finite subset of $\Omega$, $\mathcal{T}=\{0,7,24,46\}$ hours, and $N_d$ denotes the number of space-time points in $\mathcal{T} \times \Omega $ where we have observations. 
The PDE loss term $\mathcal{L}_{r}$ is defined as 
\begin{align}
    \mathcal{L}_{r} (c, D) =  
     \frac{1}{|\mathcal{P}|}\sum \limits_{(x, t) \in \mathcal{P}} |\partial_t c(x,t) - D \Delta c(x, t)|^p,
    \label{math:pdeloss}
\end{align}
where $p \in [1, \infty)$, the set $\mathcal{P}$ consists of $N_r$ 
points in $\tau \times \Omega_r$, $\tau \subset [0,T]$, and $\Omega_r \subset \Omega$ is a set of $N_p=|\mathcal{P}|$ coordinates $x \in \mathbb{R}^3$ that lie in the \emph{interior} of the domain $\Omega$. The sampling strategy to generate $\mathcal{P}$ is explained in detail in \supplref{Section}{suppl-sec:hyperparameters}.
In this work we test training with both $p=2$ and $p=1$. 
It is worth noting that boundary conditions are not included (in fact, they are often not required for inverse problems\cite{raissi_physics-informed_2019}) in the PINN loss function (\ref{math:pinn}), allowing us to sidestep making additional assumptions on the unknown boundary condition. 
The initial condition is taken to be the first image at $t=0$ and simply enters via the data loss term \eqref{math:pinn-dataloss}.
A detailed description of the network architecture and other hyperparameter settings can be found in \supplref{Section}{suppl-sec:hyperparameters}.

\subsection{The finite element approach\label{sec:methods-fem}}
Our parameter identification problem 
describes a nonlinear 
ill-posed inverse problem\cite{ito_choice_1992_2, holler_bilevel_2018, kaltenbacher_adaptive_2011}. 
We build on the numerical realization of Valnes \emph{et al.}\cite{valnes_apparent_2020} and define the 
PDE constrained optimization problem \cite{hinze2008optimization} as
\begin{align}
    \min \limits_{D, g}  \sum_{t_i \in \mathcal T} \int_\Omega (c(x, t_i; D,g) - c^d(x, t_i))^2 \, \mathrm{d}x + \frac{1}{2} \int_0^T \int_{\partial \Omega} \left(\alpha |g|^2 + \beta |\frac{\partial}{\partial t} g|^2
    + \gamma | \nabla g |^2 \right) \mathrm{d}S \mathrm{d} t,
    \label{math:fem-functional}
\end{align}
where, similar to \cite{valnes_apparent_2020}, the second term is Tikhonov regularization 
with
regularization parameters $\alpha, \beta, \gamma > 0$ and $c=c(x,t, D, g)$ solves \eqref{math:pde} with initial and boundary conditions
\begin{align}
        c(x,t) = g(x,t) \quad &\text{on } \partial \Omega \times (0, T), \label{math:pde-bc} \\
    c(x, 0) = 0 \quad &\text{in } \Omega.
    \label{math:pde-ic}
\end{align}
Since the realization is based on a reduced formulation 
using the solution operator of the partial differential equations,
the boundary condition $g: \partial \Omega \times (0,T) \to \mathbb R$ needs to be introduced as a control variable. 
A detailed description of all hyperparameter settings can be found in \supplref{Section}{suppl-sec:hyperparameters}.


\subsection{Parameterization of the diffusion coefficient \label{sec:dparam}}
Previous findings \cite{croci_uncertainty_2019, ray_analysis_2019, holter_interstitial_2017, koundal_optimal_2020} indicate that diffusion contributes at least to some degree to the distribution of tracers in the brain. 
It can thus be assumed that a vanishing diffusion coefficient is unphysical.
This assumption can be incorporated into the model by parameterizing
$D$ in terms of a trainable parameter $\delta$ as
\begin{align}
    D(\delta) = D_{\mathrm{min}} + \sigma (\delta) D_{\mathrm{max}},
    \label{math:d-bounded}
\end{align}
where $\sigma(x) = (1+\exp(-x))^{-1}$ denotes the logistic sigmoid function.
In all results reported here, we initialize with $\delta = 0$ and set $D_{\mathrm{min}} = 0.1 \, \text{mm}^2 \, \text{h}^{-1}$
and $D_{\mathrm{max}} = 1.2 \, \text{mm}^2\, \text{h}^{-1}$. 
This parameterization with a sigmoid function effectively leads to vanishing gradients $|\frac{\partial D}{\partial \delta}|$ for $|\delta| \gg 1$.
In section \ref{results:simdata-w-noise} we demonstrate that this choice of parameterization can help to avoid instabilities that occur during PINN training without parameterization, i.e. 
\begin{align}
    D = \delta.
    \label{math:d-delta}
\end{align}
The reason to introduce a $D_{\mathrm{min}} > 0$ is to avoid convergence into a bad local minimum. 
For the finite element approach, we did not observe convergence into a local minimum where $D=0$, and hence used the parameterization (\ref{math:d-delta}).

\section*{Data availability}
The datasets analyzed in the current study are available from the corresponding author upon request.

\section*{Acknowledgment}

We would like to thank Lars Magnus Valnes for insightful discussions and providing scripts for preprocessing the MRI data and meshing. We would like to thank George Karniadakis, Xuhui Meng, Khemraj Shukla and Shengze Cai from Brown University for helpful discussions about PINNs in the early stages of this work. We note and thankfully acknowledge G. Karniadakis' suggestion to switch to $\ell^1$ loss during the optimization.
The finite element computations were performed on resources provided by Sigma2 - the National Infrastructure for High Performance Computing and Data Storage in Norway.
The PINN results presented in this paper have been computed on the Experimental Infrastructure for Exploration of Exascale Computing (eX3), which is financially supported by the Research Council of Norway under contract 270053.

\section*{Author contributions}

B.Z., J.H., M.K, K.A.M. conceived the experiments. 
P.K.E. and G.R. acquired the data.
B.Z. implemented the simulators. 
B.Z. conducted the experiments and made the figures. 
All authors discussed and analyzed the results. 
B.Z., J.H., M.K., K.A.M. wrote the draft. 
All authors revised the manuscript and approved the final manuscript.

\section*{Competing interests}
The authors declare no competing interests.

\bibliography{ref2}

\clearpage

\begin{appendix}

\section{Data Generation \label{suppl-sec:data}}
\subsection{MRI Data\label{suppl-sec:mridata}}
The data under consideration in this study is based on MRI scans taken of a patient who was imaged at Oslo University Hospital in Oslo, Norway. 
The patient was diagnosed with normal pressure hydrocephalus and is referred to as "NPH1" in \cite{valnes_apparent_2020}.

The imaging protocol starts with the acquisition of a baseline MRI before 0.5 mL of a contrast agent (1 mmol/mL gadobutrol) is injected into the CSF at the spinal canal (intrathecal injection). The pulsating movement of CSF transports the tracer towards the head where it enters the brain.
The contrast agent alters the magnetic properties of tissue and fluid, and in subsequently taken MRI, enriched brain regions display changes in MR signal relative to the baseline MRI. 
From the change in signal we estimate the concentration of tracer per voxel at timepoints 0, 7, 24 and 46 hours after injection.  
Further details on the MRI acquisition and tracer concentration estimation can be found in \cite{valnes_apparent_2020}.

We next use FreeSurfer \cite{fischl_freesurfer_2012} to segment the baseline image into anatomical regions and obtain binary masks for white and gray matter. 
The human brain has many folds and represents a highly complex geometry. To limit the intrinsically high computational requirements of inverse problems, we  focus on a subregion of the white matter shown in main Fig. \ref{fig:roi}.

In the following, we describe how this data is processed further to obtain patient-specific finite element meshes to generate synthetic test data by simulation.  


\subsection{Synthetic data}
We use the surface meshes created during the brain segmentation with FreeSurfer  \cite{fischl_freesurfer_2012} to create finite element meshes of $\Omega$. 
The surface are loaded into SVMTK \cite{mardal_simulating_2022}, a Python library based on CGAL \cite{the_cgal_project_cgal_2022}, for semi-automated removal of defects and creation of high quality finite element meshes. Details on SVMTK and the mesh generation procedure can be found in \cite{mardal_simulating_2022}.

{Using FEM we then solve the PDE (\ref{math:pde}) with boundary and initial conditions (\ref{math:pde-bc}), (\ref{math:pde-ic}) 
with a diffusion coefficient $D_0 = 0.36\,$mm$^2\,$h$^{-1}$.}
This value for the diffusion coefficient of gadubutrol was estimated in \cite{valnes_apparent_2020} from diffusion tensor imaging (DTI).
In detail, we discretize (\ref{math:pde}) using the Crank-Nicolson scheme and use integration by parts to transform  (\ref{math:pde}) into a variational problem that is solved in FEniCS \cite{aln_fenics_2015} with continuous linear Lagrange elements. 
We use a high resolution mesh with $3\times 10^5$ vertices ($10^6$ cells) and small time step of $16\,$min. 
In combination with the Crank-Nicolson scheme, this minimizes effects of numerical diffusion.
For the initial condition (\ref{math:pde-ic}) we assume no tracer inside the brain at $t=0$, i.e. $c_0=0$.
The boundary condition (\ref{math:pde-bc}) is assumed to be spatially homogeneous while we let it vary in time as
\begin{align}
	g(t) = \begin{cases} 2 {t}/{T} &\text{ for } 0  \leq t \leq T/2 
	\\2 -2t/T &\text{ for } T/2  \leq t \leq T.
	\end{cases}
    \label{suppl-math:simdata-bc}
\end{align}
This choice leads to enrichment of tissue similar to what is observed experimentally over the timespan of $T=46$ hours. 
Finally, we interpolate the finite element solution $c(x,t)$ between mesh vertices and evaluate it at the center coordinates $x_{ijk}$ of the voxels $ijk$ inside the region of interest $\Omega$ and store the resulting concentration arrays $c_{ijk}$ at 0, 7, 24 and 46 hours.
With this downsampling procedure, we are then able to test the methods within the same temporal and spatial resolution as available from MRI. 

\subsection{Synthetic data with artifical noise\label{suppl-sec:noise}}
We test the susceptibility of the methods with respect to noise by adding to the data perturbations drawn randomly from the normal distribution $\mathcal{N}(0, \sigma^2 )$. We refer to the standard deviation $\sigma$
 as noise level hereafter. 
Since negative values for the concentration $c$ are nonphysical, we threshold negative values to 0, i.e. the noise-corrupted voxel values are computed as
\begin{align}
    c_{ijk} = \max \{0, c_{ijk} + \eta \} \quad \mathrm{where} \quad \eta \sim \mathcal{N}(0, \sigma^2). 
    \label{suppl-math:noise} 
\end{align}

In all the results presented in this work, we choose $\sigma = 0.05$. This corresponds to 5 \% of the maximum value of $c=1$ in the simulated measurements $c_{ijk}$ and allows to reproduce some of the  characteristic difficulties occurring when applying the PINN to the clinical data considered here.

\section{Hyperparameter settings \label{suppl-sec:hyperparameters}}
In our PINN approach, we model $c: (x, t) \rightarrow \mathbb{R}, x\in\mathbb{R}^3,~ t\in [0,T]$ by a feedforward neural network with 9 hidden layers and 64 neurons in each layer and hyperbolic tangent as activation function together with Glorot initialization \cite{glorot_understanding_2010} for all results presented here. 
We have also experimented with larger networks, different and adaptive activation functions but have not observed significant differences for a range of choices in terms of convergence rates or accuracy. 
We note that since the raw data is already a 3-D array representing grid structured data, using (physics-informed) convolutional neural networks instead of fully connected networks might yield benefits such as computational speed up. In this work, however, we decide to focus on tuning the loss function formulation and find that using a fully connected network in combination with a properly tuned loss function yields results that are consistent with the FEM approach. 

The network furthermore has an input normalization layer with fixed parameters to normalize the inputs to the range $[-1, 1]$.
To set these weights, we first compute the smallest bounding box containing all points $x=(x_1,x_2, x_3) \in \Omega$ to obtain lower and upper bounds $l_i,u_i$, $i=1, 2, 3$ such that $l_i \leq x_i \leq u_i$ for all $x\in \Omega$.
The first layer normalizes the inputs as
\begin{align}
\begin{split}
&t \hookleftarrow 2\frac{t}{T} -1, \medspace x_i  \hookleftarrow 2\frac{x_i-l_i}{u_i-l_i} -1
\end{split}
\end{align}
for $i=1,2,3$ and with $T=46\,$h (last MRI acquisition timepoint).

If not stated otherwise, we use $N_p = 10^{6}$ space-time points $(x_1,x_2,x_3,t)$ for the evaluation of the PDE loss \ref{math:pinn}). We found that this number was either sufficiently high to reach accurate recovery of the diffusion coefficient, or more sophisticated refinement techniques like residual-based adaptive refinement (RAR) \cite{lu_deepxde_2021} were needed instead of simply using more PDE points. 
The samples for the spatial coordinates $x_1,x_2,x_3$ are generated by first drawing a random voxel $i$ inside $\Omega$. 
The voxel center coordinates $(x^i_1, x^i_2, x^i_3)$ are then perturbed to obtain $x_1 = x^i_1 + dx$ where $x^i_1$ is the $x_1$-coordinate of the center of a randomly drawn voxel $i$, and similarly for $x_2$ and $x_3$. 
The perturbation $dx$ is drawn from the uniform distribution $\mathcal{U}([-0.5 \, \mathrm{mm}, 0.5 \, \mathrm{mm}]) $ and ensures that $(x_1, x_2, x_3)$ lies within the voxel $i$ (the voxels correspond to a volume of $1\, \mathrm{mm}^3$). The values for $t$ are chosen from a latin hypercube sampling strategy over the interval $[0, T]$.
We furthermore normalize the input data $c^d$ by the maximum value such that $0 \leq c^d \leq 1$. 
In both the simulation dataset and the MRI data considered here, we use four images and the same domain $\Omega$. The binary masks describing $\Omega$ consist of roughly  $0.75 \times 10^4$ voxels, i.e., the four images (at 0, 7, 24 and 46 hours) yield a total of $N_d = 3 \times 10^5$ data points. 
Due to the large number of data and PDE points, we use minibatch sampling 
of the PINN loss function (\ref{math:pinn}) and minimize it using the ADAM optimizer \cite{kingma_adam_2017}. The learning rate $\eta$ as well as potential learning rate decay schemes are an important hyperparameter, and we specify the used values in each section. 
The training set is divided into 20 batches, corresponding to $10^4$ and $5\times 10^4$ samples per minibatch in the data and PDE loss term, respectively.

Details on the implementation of the minibatch sampling strategy are presented in Section \ref{suppl-appendix:minibatching}.
It is worth noting here that our main reason to use minibatch sampling are not memory limitations. 
The graphics processing units (NVIDIA A100-SXM4) that we use to train the PINN have 80$\,$GB of memory. 
This is enough to minimize the PINN loss function (\ref{math:pinn}) with $N_d = 3 \times 10^5$ and $N_p = 10^6$ data and PDE loss points in a single batch. Our reason to use minibatch sampling is that the stochasticity of minibatch gradient descent helps to avoid local minima, 
see, e.g., Chapter 8 in \cite{Goodfellow-et-al-2016}.
In Section \ref{suppl-sec:simdata-clean} below we perform a systematic study using different minibatch sizes and find that smaller batch sizes are preferable in our setting since they yield more accurate recovery of the diffusion coefficient (for a fixed number of epochs).

As for the finite element approach, we discretize (\ref{math:pde}) in time using the Crank-Nicolson scheme and 48 time steps. We then formulate the PDE problem as variational problem and solve it in FEniCS \cite{aln_fenics_2015} using the finite element method. 
To limit the compute times required, we use a time step size of $1\,$h and continuous linear Lagrange elements. 
We further use linear interpolation of the data as a starting guess for the boundary control $g$, 
\begin{align}
	g(x,t) = c^d(x, t_i) + \frac{c^d(x, t_{i+1})-c^d(x, t_i)}{t_{i+1} -t_i} (t- t_i) \label{suppl-math:opt-control-starting-guess}
\end{align}
for $t_i \in \mathcal{T} = \{0,7,24,46\}\,$h and $t_i \leq t \leq t_{i+1}$.
We then use dolfin-adjoint \cite{mitusch_dolfin-adjoint_2019} to compute gradients of the functional (\ref{math:fem-functional}) with respect to $D$ and $g$ and optimize using the L-BFGS method. 

In terms of degrees of freedom (optimization parameters), these settings result in 33665 weights in the neural network. For the finite element approach, the degrees of freedom depend on the number of vertices on the boundary of the mesh since the control is the boundary condition $g$. Our mesh for $\Omega$ has 33398 cells on the boundary. 
For 48 time steps, this yields $48\times 33398 = 1.6 \times 10^6 $ degrees of freedom.

\section{Validation on synthetic data \label{suppl-sec:simdata-clean}}

We verify the implementation of the two approaches by considering synthetic data without noise, cf. Fig. \ref{fig:068datasim}, in the white matter subregion $\Omega$ depicted in Fig. \ref{fig:roi}.

For the PINN approach we test different minibatch sizes for three different optimization schemes
using the ADAM optimizer: (i) fixed learning rate $10^{-3}$ and $p=2$, (ii) fixed learning rate $10^{-3}$ while we switch from $p=2$ to $p=1$ after half the epochs and (iii) using initial learning rate $10^{-3}$ that decays exponentially during training to $10^{-4}$ and $p=2$. 
Table \ref{suppl-table:minibatch-pinn-clean} tabulates the relative error between the learned diffusion coefficient and the ground truth $D_0$ for a wide range of parameters. We find that (a) in general smaller batch sizes result in more accurate results and (b) the results are both most stable and accurate when using exponentially decaying learning rate. 
Notably, the PINN recovers the ground truth diffusion coefficient $D_0$ to up to 1$\,$\% accuracy when using the learning rate decay optimization scheme.
These result are in line with \cite{kadeethum_physics-informed_2020} where increased accuracy in parameter recovery was observed for smaller batch sizes. However, there are also settings where full batch optimization with L-BFGS improves PINN performance in parameter identification problems \cite{mathews_uncovering_2021}.

For the finite element approach, Table \ref{suppl-tab:opt-control-clean-32} presents the accuracy of the recovered diffusion coefficient.
According to the theoretical results, decreasing regularization parameters leads to higher accuracy but less well conditioned optimization problems. This is in line with the results presented in Table \ref{suppl-tab:opt-control-clean-32}. The finite element approach with appropriate regularization parameters and the PINN approach yield comparably accurate results.

\subsection{Computational effort\label{suppl-sec:compute-effort}}

We here list the computing times to estimate the diffusion coefficient from noisy simulation data with our implementation of the FEM and PINN approaches presented in the main text.

Our implementation of the FEM approach using dolfin-adjoint \cite{mitusch_dolfin-adjoint_2019} using a mesh with 91849 cells and 23307 vertices (whereof 33398 and 16693 are on the boundary, respectively) requires around 45-48 hours computing time for the 1,000 iterations until convergence as shown in Supplementary Fig. \ref{suppl-fig:cleansimdataR32convergence} using a single Intel Xeon Gold CPU.

As for the PINN approach, we terminated the optimization after 2,000 epochs of ADAM with data and PDE batch sizes of $1.5 \times 10^4$ and $5 \times 10^4$.
With our implementation in PyTorch\cite{paszke2019pytorch2} this takes about 6 hours on a NVIDIA A100-SXM4.

We note that neither of these implementations have been optimized to reduce the compute times. 

\subsection{PINN solution of the synthetic testcase}

 \begin{table}[H]
\centering
\caption{
Average rel. error $|D_{\mathrm{pinn}}-D_{0}|/D_{0}$ in \% after $2 \times 10^{4}$ epochs training on synthetic data without noise, with Algorithm 1. We average over 5 runs, numbers in brackets are standard deviation.
\label{suppl-table:minibatch-pinn-clean}}
\begin{tabularx}{\columnwidth}{ @{}|X|X|XXX|@{} }
\hline
Optimization scheme & \diagbox{$n_r$}{$n_d$}  & $10^4$  &  $5 \times 10^4$ & $ 10^5$  \\
 \hline
 \multirow{5}{*}{\makecell{ADAM \\ lr = 1e-3 \\ $p=2$} } 
 & 10000  & 2 (0) & 2 (1)   & 3 (1)   \\
 & 33334  & 4 (0) & 12 (6)  & 9 (1)   \\
 & 50000  & 7 (1) & 5 (0)   & 2 (0)   \\
 & 100000 & 7 (1) & 55 (18) & 59 (18) \\
 & 166667 & 8 (1) & 24 (17) & 50 (23) \\
 & 333334 & 8 (1) & 38 (4)  & 50 (16) \\
\hline
 \multirow{5}{*}{\makecell{ADAM \\ lr = 1e-3 \\ $p$=$2 \rightarrow$ $p$=1}}
 & 10000  & 2 (0) & 1 (0)   & 2 (0)   \\
 & 33334  & 2 (0) & 2 (0)   & 2 (0)   \\
 & 50000  & 2 (1) & 2 (0)   & 2 (1)   \\
 & 100000 & 2 (1) & 58 (27) & 50 (25) \\
 & 166667 & 2 (1) & 3 (2)   & 68 (4)  \\
 & 333334 & 2 (0) & 0 (0)   & 62 (6)    \\
\hline
 \multirow{5}{*}{\makecell{ADAM \\ exp lr decay \\ 1e-3 $\rightarrow$ 1e-4 \\ $p=2$}}
  & 10000  & 1 (0) & 1 (0)  & 1 (0)   \\
 & 33334  & 1 (0) & 1 (1)  & 1 (0)   \\
 & 50000  & 1 (1) & 1 (0)  & 1 (0)   \\
 & 100000 & 1 (0) & 10 (6) & 72 (0)  \\
 & 166667 & 1 (1) & 4 (2)  & 23 (26) \\
 & 333334 & 1 (0) & 4 (4)  & 59 (23) \\
\hline
\end{tabularx}
\end{table}

\subsection{Finite element solution of the synthetic testcases\label{suppl-sec:fem-synthetic}}

\begin{table}[H]
\centering
\caption{Rel. error $|D-D_0|/D_0$ for the FEM approach (\ref{math:fem-functional}), different regularization parameters, 3 measurement points, clean data, 1,000 iterations. Convergence of the optimization is demonstrated in Fig. \ref{suppl-fig:cleansimdataR32convergence}.
\label{suppl-tab:opt-control-clean-32}}
\begin{tabular}{ll|llll}
\hline
$\alpha$ & \diagbox[]{$\beta$}{$\gamma$} & 0.0  & 0.01 & 1.0  \\
\hline
  & 0.001 & 9  & 8 & 261 \\
$10^{-6}$  & 0.01  & 1  & 5 & 261 \\
  & 0.1   & 11 & 10 & 261  \\
\hline
\end{tabular}
\end{table}

\begin{table}[H]
\centering
\caption{Rel. error $|D-D_0|/D_0$ in \%  for the finite element approach (\ref{math:fem-functional}), different regularization parameters.
Failure of the algorithm is indicated by the symbol "x". Convergence plots for the optimization are given in Fig. \ref{suppl-fig:noisyimdataR32convergence}.
\label{suppl-tab:opt-control-w-noise}}
\begin{tabular}{ll|ll}
\hline
$\alpha$ & \diagbox[]{$\beta$}{$\gamma$} & 0.0 & 0.01 \\
\hline
  & 0.0  & 43  & 10   \\
$10^{-6}$  & 0.01 & 8   & x    \\
  & 0.1  & 4   & x    \\
  \hline
 & 0.0  & 44  & x    \\
$10^{-4}$ & 0.01 & 5   & 13   \\
 & 0.1  & 6   & 12  \\
\hline
\end{tabular}
\end{table}

\begin{figure}[H]
\centering
     \begin{subfigure}[b]{0.46\textwidth}
         \centering
         \includegraphics[width=\textwidth]{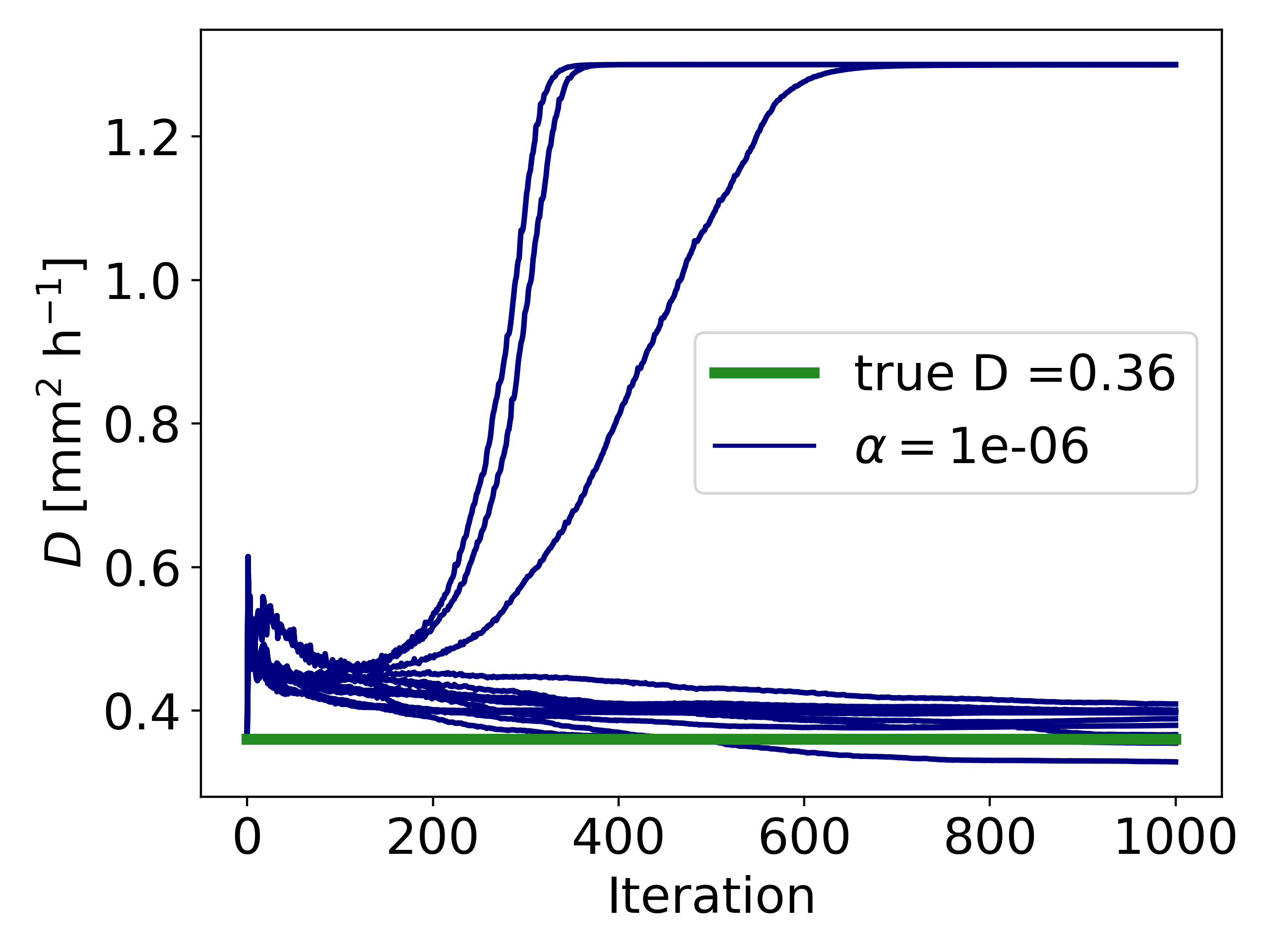}
         \caption{}
         \label{suppl-fig:cleansimdataR32convergence}
     \end{subfigure}
     \begin{subfigure}[b]{0.46\textwidth}
         \centering
         \includegraphics[width=\textwidth]{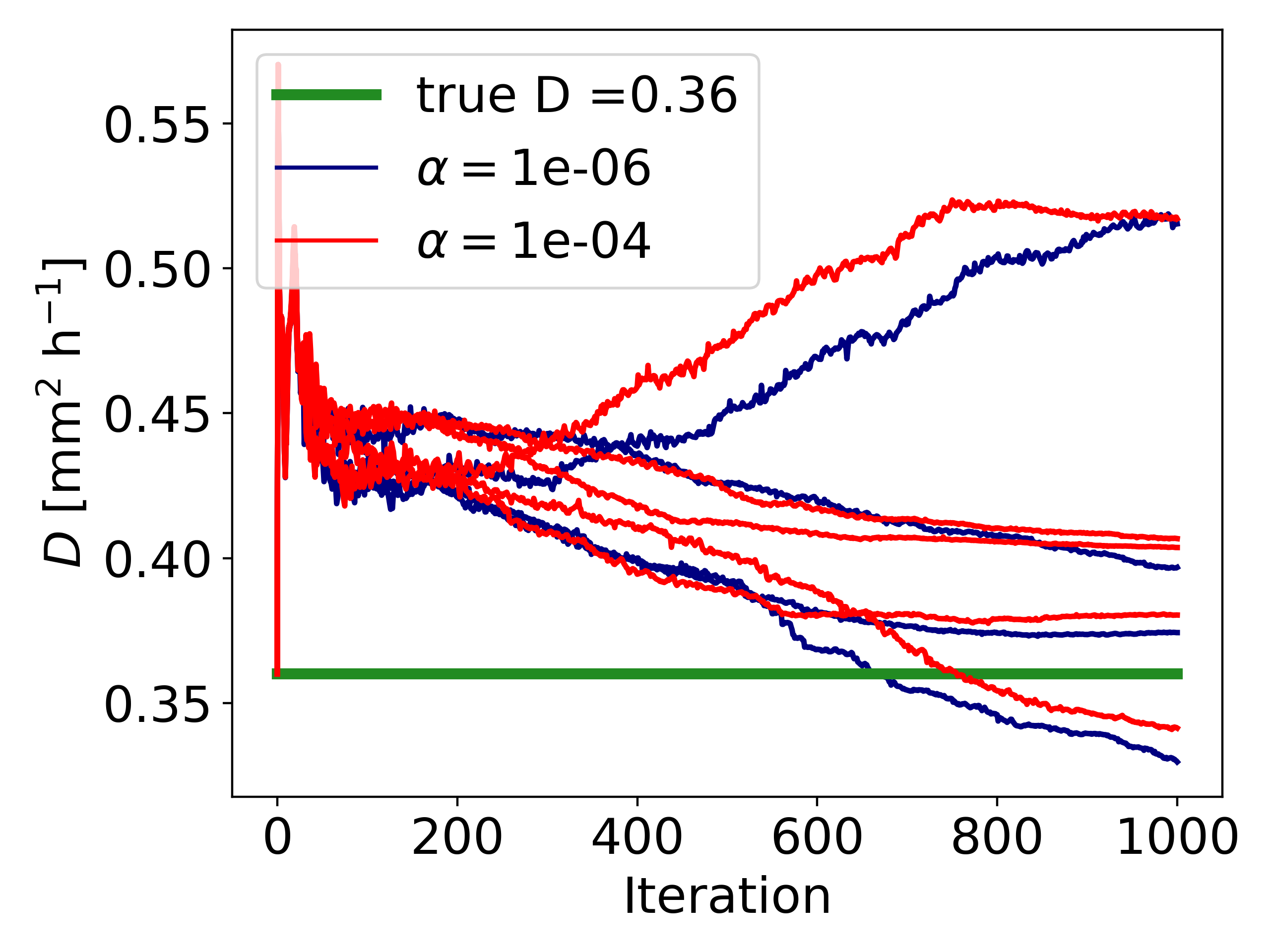}
         \caption{\label{suppl-fig:noisyimdataR32convergence}}
     \end{subfigure}
     \caption{
        Convergence plots for the FEM regularization parameters presented in Tables \ref{suppl-tab:opt-control-clean-32} (synthetic data)
         and \ref{suppl-tab:opt-control-w-noise} (synthetic data with noise).
         }
\end{figure}

\section{Details on the PINN training procedures \label{suppl-appendix:minibatching}}

\subsection{Minibatch sampling strategy}

\begin{algorithm}[H]
 \caption{Minibatch training
 \label{suppl-alg:minibatching}}
 \begin{algorithmic}[1]
 \renewcommand{\algorithmicrequire}{\textbf{Input:}}
 \renewcommand{\algorithmicensure}{\textbf{Output:}}
 \REQUIRE neural network $c$ with parameters $\theta$,
 data minibatch size $n_d$,
 PDE minibatch size $n_r$, 
 RAR checkpoints $\{i_1, \dots, i_n \}$,
 epochs, learning rate $\eta$, initial guess $D$ for the diffusion coefficient,
   input-data pairs $\lbrace ( \mathbf{x}_k^d, c_k^d)  \in \mathbb{R}^4 \times \mathbb{R} \rbrace_{1 \leq k \leq N_d}$, 
   PDE space-time points \begin{math} \mathcal{P}=\lbrace \mathbf{x}_{k} \in \mathbb{R}^4 \rbrace_{1 \leq k \leq N_r} \end{math}
  
  \STATE compute number of data batches $b_d = $ceil$ (N_d / n_d)$
\STATE compute number of PDE batches $b_r = $ ceil$ (N_r / n_r)$
\STATE Set $b = $max$(b_d, b_r)$
  \FOR {$i$ in range(epochs)}
  \IF{i $\in$ RAR checkpoints}
  \STATE add points to $\mathcal{P}$ with either procedure \ref{suppl-alg:rar} or \ref{suppl-alg:rae}
  \ENDIF 
  \STATE randomly split $\lbrace ( \mathbf{x}_k^d, c_k^d) \rbrace$ into subsets $\mathcal{D}_{1 \leq j \leq b_d}$
\STATE randomly split $\lbrace \mathbf{x}_r^k \rbrace$ into subsets $\mathcal{R}_{1 \leq j \leq b_r}$
  \STATE \# Iterate over all minibatches
    \FOR {$j$ in range$(b)$}
    \STATE \# Start from beginning should you reach the last subset in $\mathcal{D}_{b_k}$ or $\mathcal{R}_{b_r}$, respectively (Happens if $b_r \neq b_d$):
    \STATE Set $j_d = j \, \mathrm{mod} \, b_d$, $j_r = j \, \mathrm{mod} \, b_r$
    \STATE \# Compute losses on subsets
    \STATE $\mathcal{L}$ = $\frac{1}{|\mathcal{R}_{j_d}|} \sum \limits_{\mathbf{x}^d, c^d \in \mathcal{D}_{j_d}} \left(c(\mathbf{x}^d) - c^d \right)^2$
    \STATE $\mathcal{L}$ += $\frac{1}{|\mathcal{R}_{j_r}|} \sum \limits_{\mathbf{x} \in \mathcal{R}_{j_r}} \left|\partial_t c(\mathbf{x}) - D \Delta c(\mathbf{x})\right|^p$
    \STATE \# update parameters $\theta$
    \STATE $\theta$ -= $\eta \nabla_{\theta} \mathcal{L}$
    \STATE \# update diffusion coefficient $D$
    \STATE $D$ -= $\eta \nabla_{D} \mathcal{L}$
    \ENDFOR
  \ENDFOR
 \end{algorithmic} 
 \end{algorithm}

\subsection{Residual based refinement 
	\label{suppl-app:RAR-RAE}}

\begin{algorithm}[H]
 \caption{Refinement step with the RAR procedure as in "Procedure 2.2" in \cite{lu_deepxde_2021} adapted to the nomenclature in our work.
 \label{suppl-alg:rar}}
 \begin{algorithmic}[1]
  \renewcommand{\algorithmicrequire}{\textbf{Input:}}
 \renewcommand{\algorithmicensure}{\textbf{Output:}}
 \REQUIRE The set of {$N_r$} PDE points $\mathcal{P}$, PDE residual $r(x,t)$, number $m$ of points to add per refinement step, number $n$ of points to test the residual 
\ENSURE refined set of PDE points $\mathcal{P}$
\STATE Compute the absolute value of the PDE residual $|r(x,t)|$ at $n$ random samples $\mathcal{S} = \{(x_1, t_1), \dots, (x_n, t_n)\}$ from $\Omega_r \times {\tau}$
\STATE Sort $\mathcal{S}$ by decreasing residual $|r(x,t)|$ and keep only the first $m$ points in $\mathcal{S}_m$
\RETURN The {refined} set of $N_r+m$ points $\mathcal{P} \cup \mathcal{S}_m$ 
  \end{algorithmic} 
 \end{algorithm}

\begin{figure}[H]
\centerline{\includegraphics[width=\textwidth]{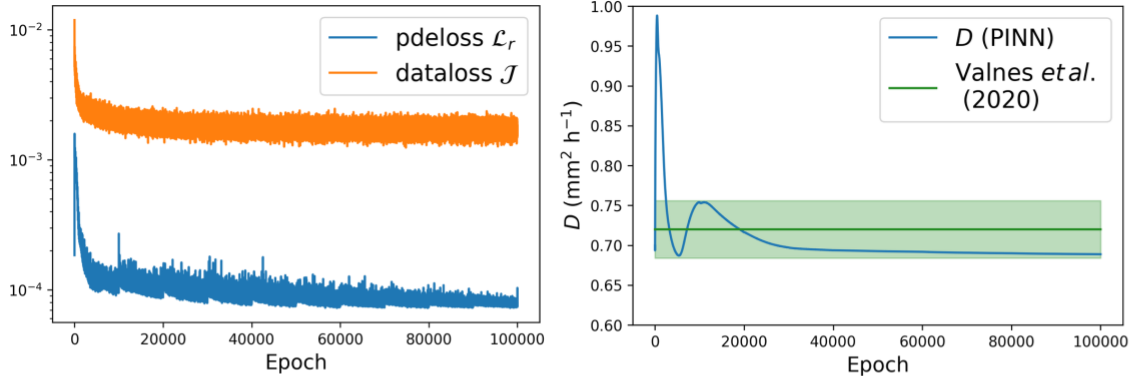}}
\caption{(Left) Convergence plots for the PINN losses trained with RAR on MRI data. (Right) Diffusion coefficient during PINN training.
Exponential learning rate decay from $10^{-4}$ to $10^{-5}$ with RAR and $p=1$.
\label{suppl-fig:cleansimdataR32convergencepinn}}
\end{figure}

\begin{algorithm}[H]
 \caption{Refinement step with the RAE procedure, a modification of the RAR procedure as described under "Procedure 2.2" in \cite{lu_deepxde_2021}.
 \label{suppl-alg:rae}}
 \begin{algorithmic}[1]
  \renewcommand{\algorithmicrequire}{\textbf{Input:}}
 \renewcommand{\algorithmicensure}{\textbf{Output:}}
 \REQUIRE The set of {$N_r$} PDE points $\mathcal{P}$, PDE residual $r(x,t)$, number $m$ of points to add per refinement step, number $n$ of points to test the residual 
\ENSURE refined PDE points $\mathcal{P}$
\STATE Compute the absolute value of the PDE residual $|r(x,t)|$ at $n$ random samples $\mathcal{S} = \{(x_1, t_1), \dots, (x_n, t_n)\}$ from $\Omega_P \times \tau$
\STATE Sort $\mathcal{S}$ by decreasing residual $|r(x,t)|$ and keep only the first $m$ points in $\mathcal{S}_m$
\STATE Compute the PDE residual $|r(x,t)|$ at the points in $\mathcal{R}$ 
\STATE Sort $\mathcal{R}$ by increasing residual $|r_{\mathcal{R}}(x,t)|$ and keep only the first $N_r-m$ points in $\mathcal{R}_{N_r-m}$
\RETURN The set of $N_r$ refined points $\mathcal{R}_{N_r-m} \cup \mathcal{S}_m$ 
  \end{algorithmic} 
 \end{algorithm}
 \end{appendix}

\end{document}